\documentclass{article}

\usepackage{arxiv}
\usepackage[table,xcdraw]{xcolor}
\usepackage{amsmath,amsthm,amssymb,amscd}
\usepackage[all,cmtip]{xy}
\usepackage{diagxy}
\usepackage[utf8]{inputenc} 
\usepackage[T1]{fontenc}    
\usepackage{hyperref}       
\usepackage{url}            
\usepackage{booktabs}       
\usepackage{amsfonts}       
\usepackage{nicefrac}       
\usepackage{microtype}      
\usepackage{lipsum}
\usepackage{mathrsfs}
\usepackage{fancyhdr}
\usepackage{color}
\usepackage{subfig}
\usepackage{graphicx}
\usepackage{textcomp}
\usepackage{ytableau}
\usepackage{tikz}
\usetikzlibrary{calc,intersections,through,backgrounds,patterns}
\usepackage{algorithm}
\usepackage{algpseudocode}
\newtheorem{exam.}{Example}
\newtheorem{def.}{Definition}

\newtheorem{prop.}{Proposition}
\newtheorem{coro}{Corollary}
\newtheorem{theorem}{Theorem}
\newtheorem*{remark}{Remark}
\newtheorem*{notation}{Notation}

\newcommand{\Obj}{{\rm{Ob}}}
\newcommand{\id}{{\rm{id}}}

\newcommand{\bbZ}{\mathbb{Z}}
\newcommand{\bbR}{\mathbb{R}}
\newcommand{\bbN}{\mathbb{N}}

\newcommand{\Hom}{{\rm{Hom}}}

\newcommand{\im}{{\rm{im}}}

\newcommand{\bfe}{{\mathbf e}}

\newcommand{\bfx}{{\mathbf x}}

\newcommand{\conv}{{\rm conv}}

\newcommand{\cl}{{\rm cl}}

\title{Locating topological structures in digital images via local homology}

\author{
  Chuan-Shen Hu \\ 
  School of Physical and Mathematical Sciences\\
  Nanyang Technological University\\
  50 Nanyang Avenue 639798, Singapore\\
  \texttt{chuanshen.hu@ntu.edu.sg}\\
  \texttt{peterbill26@hotmail.com}
}

\begin{document}
\maketitle

\begin{abstract}
Topological data analysis (TDA) is a rising branch in modern applied mathematics. It extracts topological structures as features of a given space and uses these features to analyze digital data. Persistent homology, one of the central tools in TDA, defines persistence barcodes to measure the changes in local topologies among deformations of topological spaces. Although local spatial changes characterize barcodes, it is hard to detect the locations of corresponding structures of barcodes due to computational limitations. The paper provides an efficient and concise way to divide the underlying space and applies the local homology of the divided system to approximate the locations of local holes in the based space. We also demonstrate this local homology framework on digital images.
\end{abstract}


\section{Introduction}

Homology is an algebraic description of topological spaces and has become one of the foundations of modern geometry and topology. It uses algebra to detect genera in topological spaces, such as loops and high-dimensional voids, and to classify the topological types and shapes of manifolds. In addition to its importance in pure mathematics, over the past two decades or so, data scientists have noticed the benefits and potential of homology in numerical data and raised a new field called topological data analysis (TDA)~\cite{zomorodian2004computing, carlsson2009topology,ghrist2008barcodes,carlsson2010multiparameter, garside2021event}.

Persistent homology plays a central role in TDA, which transforms a sequence of topological spaces linked by continuous functions into a homology chain. By checking the birth and death of elements in the chain, one can understand which homological generator can have a longer lifespan and shows its importance in the continuous process~\cite{zomorodian2004computing}.  Persistent homology and related techniques have been applied in many data science tasks, such as bioinformatics~\cite{riihimaki2020topological,nielson2015topological,hu2021toporesnet}, molecular analysis~\cite{xia2014persistent,gong2022persistent,anand2022topological,wee2022persistent}, image processing~\cite{chung2022multi,chung2018topological, edelsbrunner2013persistent,de2023value,onuchin2023individual}, and material science~\cite{chung2018topological}.

Persistent barcodes (Section \ref{Subsec. Persistent Homology}) record the lifespans of connected components, loops, and voids. Many applications use persistence barcodes and related statistical features as machine learning features~\cite{bubenik2015statistical,adams2017persistence,chung2022persistence}. Although persistent homology and persistence barcode has shown their potential in many real applications, it still has some limitations. One is it can only capture the global information of how connected components and holes behave during geometric deformation, while the local merging relations are usually omitted. This information is theoretically present in the definition of persistent homology and persistent barcodes. However, for computational efficiency, hole representations (e.g., $q$-circular representation of $q$-holes) or positions are often buried in the Gaussian elimination of the matrices in the computation of persistence barcodes. 

Recently, some scholars noticed the importance of local information on persistent homology and proposed some interesting works on the local behavior of persistent homology~\cite{vandaele2019local, stolz2019global}.  For example, Vandaele \textit{et al.}~\cite{vandaele2019local} investigated the local Vietoris-Rips complexes of the point cloud and applied the local Betti pairs to form a global descriptor of the point cloud. This descriptor can be viewed as a heatmap of the whole space. Regions with higher heat values usually mean they have more significant topological/geometric information, such as higher local branch numbers or loops. Also, Stolz described in her doctoral dissertation~\cite{stolz2019global} how to apply the Mayer--Vietoris sequence to compute the local Vietoris--Rips complex linked from a data point.

On the other hand, some of the research also aims to detect the locations of loop or hole structures in the topological space. For example, Akai \textit{et al.}~\cite{9575312} generate persistence barcodes of the Vietoris--Rips complex as inputs of a neural network model and apply them for the ego-vehicle localization application. Similarly, Keros \textit{et al.}~\cite{keros2022dist2cycle} train on a Hodge Laplacian-based graph neural network to detect the nearest optimal homology as a location representation of homologies. Furthermore, Xu \textit{et al.}~\cite{xu2019finding} apply the distance measurement (DTM) function~\cite{chazal2017robust} to enhance the robustness of Vietoris--Rips complex construction, and apply persistence and distance information to detect holes and voids in point cloud data.

However, while image structures are more regular than point clouds, making it easier to compare local-global attributes, most current methods are designed for point cloud data. Theoretical assurance methods for localized hole detection are still limited. The paper provides a theoretically guaranteed framework for hole position detection in arbitrary topological spaces, demonstrated on digital images. 

This paper is an extension of our previous work presented as a workshop paper at CVPR 2021 (2021 Conference on Computer Vision and Pattern Recognition)~\cite{hu2021sheaf}. The work introduces the concept of cellular sheaves and connects them to persistent homology. In this work, we define the \textit{local merging number} and consider its geometric significance in \( 0 \)-dimensional objects. This paper extends the framework to focus on merging relations in \( 1 \)-dimensional structures. In addition to theoretical promotion, we have a preliminary demonstration of images. It shows that the 1-dimensional merging relations can estimate the position of holes in the space, which provides a way to analyze the local topological characteristics. 

\subsection{Organization}

The organization of the paper is as follows. Section \ref{Sec. Persistent Homology and Barcodes} quickly recaps the homology, Betti numbers, persistent homology, and barcodes. We present the main results in Section \ref{Sec. Our Approaches} and separate the section into two parts. Section \ref{Subsec. Persistent Homology of Local Systems} introduces how we divide the ambient space by a local region and apply the divided system to compute its persistent homology. We also interpret the geometric meaning of the computed barcodes and explain how they detect the cycle locations.   We also compare the proposed framework with previous methods in Section \ref{Subsec. Discussion}. Section \ref{Sec. Demonstration} shows how to adapt the theory developed in Section \ref{Sec. Our Approaches} on digital images and demonstrates the proposed locating method. Finally, we discuss future directions and summarize the paper in Section \ref{Sec. Conclusion}.

\section{Persistent Homology and Barcodes}
\label{Sec. Persistent Homology and Barcodes}
We briefly introduce the standard notions and terminologies of singular homology, including its functoriality, Betti numbers, and geometric meanings in Section \ref{Subsec. Homology}. Section \ref{Subsec. Persistent Homology} focuses on persistent homology and barcodes. We will also show in this section typical ways for building filtrations, especially the construction relying on the thresholding technique, which is the foundation of the paper.
\subsection{Homology}
\label{Subsec. Homology}
This section briefly recalls the singular homology and related properties of topological spaces. One can find these materials in several classic textbooks on algebraic topology~\cite{Hatcher:478079,Vick,munkres2018elements,Greenberg}. We start the section with the following definitions.

\begin{def.}
For any non-negative integer $q$, we define the \textbf{geometric }$\textbf{q}$\textbf{-simplex}, denoted by $\Delta_q$, as the convex hull of the standard basis $\{ \bfe_0, \bfe_1, ..., \bfe_q \}$ for the $(q+1)$-dimensional Euclidean space $R^{q+1}$. That is, 
\begin{equation*}
    \Delta_q = \conv(\bfe_0, \bfe_1, ..., \bfe_q) = \left \{ t_0\bfe_0 + t_1 \bfe_1 + \cdots + t_q\bfe_q :  t_i \in [0,1] \ \text{and} \ \sum_{i = 0}^q t_i = 1 \right \}.
\end{equation*}
\end{def.}

For any $(q+1)$ points $\bfx_0, ..., \bfx_q$ in $\bbR^{n}$ we can define the \textit{affine map} $[\bfx_0, ..., \bfx_q] : \Delta_q \rightarrow \bbR^{n}$ by
\begin{equation}
\label{Def. Affine map}
    (t_0, t_1, ..., t_q) \longmapsto t_0\bfx_0 + \cdots + t_q \bfx_q.
\end{equation}
Then $[\bfx_0, ..., \bfx_q]$ is a continuous map. A continuous function from $\Delta_q$ to a topological space $X$ is called a \textit{singular} $q$\textit{-simplex} in $X$. In particular, any affine map $[\bfx_0, ..., \bfx_q] : \Delta_q \rightarrow \bbR^{n}$ is a singular $q$-simplex in $\bbR^n$. For $q \in \bbZ_{\geq 0}$ and $i \in \{ 0, 1, ..., q+1 \}$ we define $f_{q+1}^i = [\bfe_0, ..., \widehat{\bfe_{i}}, ..., \bfe_{q+1}]$. In other words, $f_{q+1}^i$ is a singular $(q-1)$-simplex in $\bbR^{q+1}$. One can see that the image of $f_{q+1}^i$ is actually the convex hull of the set $\{ \bfe_0, ..., \widehat{\bfe_{i}}, ..., \bfe_{q+1} \}$, which is the $i$-th $(q-1)$-face of the geometric simplex $\Delta_q$~\cite{Greenberg}. 

\begin{def.}
Let $X$ be a topological space, $R$ a commutative ring with identity, and $q$ a non-negative integer. We define $S_q(X;R)$ as the free $R$-module generated by all continuous maps $\sigma : \Delta_q \rightarrow X$. For convenience, we usually define $S_q(X;R) = 0$ for $q < 0$.     
\end{def.}

The singular simplices give us a way to express geometric simplices in arbitrary topological spaces. In Euclidean spaces, one can explore the faces as boundaries of geometric simplices by using convex analysis, while it is not applicable in general spaces. In algebraic topology, we use the following boundary maps to read the boundary data of singular simplices.

\begin{def.}
Let $X$ be a topological space, $R$ a commutative ring with identity, and $q \in \bbN$ a positive integer. The $\textbf{q}$\textbf{-boundary map} is the function $\partial_q : S_q(X;R) \rightarrow S_{q-1}(X;R)$ that extends by the mapping
\begin{equation*}
\sigma \longmapsto \sum_{i = 0}^{q} (-1)^i \cdot \sigma \circ f_q^i    
\end{equation*}
for all continuous $\sigma : \Delta_q \rightarrow X$. Note that $\partial_q$ is well-defined since each $\sigma \circ f_q^i$ is a singular $(q-1)$-simplex in $X$.
\end{def.}

Because $S_q(X;R)$ is defined as the zero space for $q < 0$, we also define $\partial_q$ as the zero maps for $q \leq 0$. The following proposition is the foundation of homology theory.

\begin{prop.}
[\cite{Greenberg}, (9.2)]
\label{Prop. Boundary Operator Property}
Let $X, R, q$ and $\partial_q$ be defined as above. Then $\partial_{q-1} \circ \partial_q = 0$.
\end{prop.}

The equation $\partial_{q-1} \circ \partial_q = 0$ shows that $\im(\partial_{q}) \subseteq \ker(\partial_{q-1})$ for every $q \in \bbZ$, and hence we can define the $q$-th singular homology of $X$ as the $R$-module
\begin{equation*}
    H_q(X;R) = \frac{\ker(\partial_{q})}{\im(\partial_{q + 1})}.
\end{equation*}
\begin{notation}
[\cite{Greenberg,Vick,munkres2018elements,EdelsbrunnerHarerbook2010}]
To simply the notations, for a topological space $X$ and $q \geq 0$, we use $Z_q(X;R)$ and $B_q(X;R)$ to denote the modules $\ker(\partial_q)$ and $\im(\partial_{q+1})$. That is, 
\begin{equation}
\label{Eq. Z/B Notations of H}
    Z_q(X;R) = \ker(\partial_q), B_q(X;R) = \im(\partial_{q+1}), \text{ and } H_q(X;R) = \frac{Z_q(X;R)}{B_q(X;R)}.
\end{equation}
Chains in $Z_q(X;R)$ and $B_q(X;R)$ are called the $\textbf{q}$\textbf{-cycles} and $\textbf{q}$\textbf{-boundaries} of $X$.  
\end{notation}
Except for sending each topological space $X$ to an $R$-module $H_q(X;R)$,  for every continuous map $f : X \rightarrow Y$ and $q \in \bbZ_{\geq 0}$ we can define an $R$-module homomorphism $S_q(f;R) : S_q(X;R) \rightarrow S_q(Y;R)$ that extends the mapping
\begin{equation*}
    \sigma \longmapsto f \circ \sigma
\end{equation*}
for all singular $q$-simplices $\sigma : \Delta_q \rightarrow X$. Note that the mapping is well-defined since $f \circ \sigma$ is also a continuous map from $\Delta_q$ to $Y$. This observation leads to the following proposition.
\begin{prop.}
[\cite{Greenberg}]
Let $\mathfrak{Top}$ and $\mathfrak{Mod}_R$ be the categories of topological spaces and $R$-modules. For each $q \in \bbZ_{\geq 0}$, the assignments $X \in \Obj(\mathfrak{Top}) \mapsto S_q(X;R)$ and $f \in \Hom_{\mathfrak{Top}}(X,Y) \mapsto S_q(f;R)$ form a functor from $\mathfrak{Top}$ to $\mathfrak{Mod}_R$.
\end{prop.}
In fact, for a continuous map $f: X \rightarrow Y$, one can prove that the rectangles in the ladder
\begin{equation*}
\xymatrix@+1.5em{
                & \cdots
                \ar[r]^{}
				& S_{q+1}(X;R)
				\ar[r]^{\partial_{q+1}(X)}
				\ar[d]^{S_{q+1}(f;R)}
                & S_{q}(X;R) 
                \ar[r]^{\partial_{q}(X)}
                \ar[d]^{S_{q}(f;R)}
                & S_{q-1}(X;R)
                \ar[r]^{}
                \ar[d]^{S_{q-1}(f;R)}
                &
                \cdots
                \\
                & \cdots
                \ar[r]^{}
        		& S_{q+1}(Y;R)
        		\ar[r]^{\partial_{q+1}(Y)}
				& S_{q}(Y;R)
				\ar[r]^{\partial_{q}(Y)}
				& S_{q-1}(Y;R)
				\ar[r]^{}
				&
				\cdots
}
\end{equation*}
of $R$-modules and $R$-module homomorphisms are commutative. Therefore, for every $q$, this ladder induces an $R$-module homomorphism
\begin{equation*}
    H_q(f;R): H_q(X;R) \longrightarrow H_q(Y;R)
\end{equation*}
that sends each equivalence class $[c]$ in $H_q(X;R)$ to the class $[S_q(f;R)(c)]$ in $H_q(Y;R)$. Furthermore, we can see that the assignment $H_q(\cdot;R)$ of topological spaces and continuous maps also forms a functor from $\mathfrak{Top}$ to $\mathfrak{Mod}_R$:
\begin{prop.}[\cite{Greenberg}]
\label{Proposition: functoriality of singular homology}
Let $\mathfrak{Top}$ and $\mathfrak{Mod}_R$ be the categories of topological spaces and $R$-modules. For each $q \in \bbZ_{\geq 0}$, the assignments $X \in \Obj(\mathfrak{Top}) \mapsto H_q(X;R)$ and $f \in \Hom_{\mathfrak{Top}}(X,Y) \mapsto H_q(f;R)$ form a functor from $\mathfrak{Top}$ to $\mathfrak{Mod}_R$.
\end{prop.}
An important purpose of developing singular homology is to detect holes in a topological space in any dimension. This property of singular homology is sometimes called the \textit{Poincar\'e  lemma} of singular homology. We state this lemma as follows.

\begin{prop.}
[Corollary (15.5),~\cite{Greenberg}]
Let $n \geq 1$ be a positive integer, and let
\begin{equation*}
    S^n = \{ (x_1, x_2, ..., x_{n+1}) \in \bbR^{n+1} : x_1^2 + \cdots + x_{n+1}^2 = 1 \}
\end{equation*}
be the $n$-sphere in $\bbR^{n+1}$. Then, for every commutative ring $R$ with identity and a non-negative integer $q \geq 0$, we have
\begin{equation}
H_q(S^n;R) \simeq 
    \begin{cases}
 	R &\null\text{ if } q = n \text{ or } q = 0, \\
 	0  &\null\text{ otherwise.}
 	\end{cases}
\end{equation}

In particular, for every topological space $X$, we have $H_0(X;R) \simeq R^m$, where $m$ is the number of path-connected components of $X$, and each path-connected component of $X$ can be represented by a constant function from $[0,1]$ to $X$. 
\end{prop.}

The Poincar\'e lemma provides us with a reliable measurement to detect the number of $q$-dimensional holes in a topological space.  This number is called the $q$-th \textit{Betti number}.

\begin{def.}
[\cite{Greenberg}]
Let $R$ be a PID. For any topological space $X$ and integer $q \geq 0$, we define the $\textbf{q}$-${\rm \textbf{th}}$ \textbf{Betti number} $\beta_q = \beta_q(X)$ of $X$ to be the rank of the $R$-module $H_q(X;R)$. In particular, when $R = F$ is a field, we have $\beta_q = \dim_F \ H_q(X;R)$.  
\end{def.}

In applications, we often set $R$ as the binary field $\bbZ_2 = \bbZ/2\bbZ$ and simplify the notation $H_q(X;\bbZ_2)$ to $H_q(X)$. In the paper, we will focus on homology over $\bbZ_2$ and the singular homology of (binary) images (see Section \ref{Sec. Our Approaches}).


\subsection{Prescient Homology}
\label{Subsec. Persistent Homology}

\begin{figure*}
    \begin{center}
    \subfloat[$f = f_0$]{\includegraphics[width=25mm]{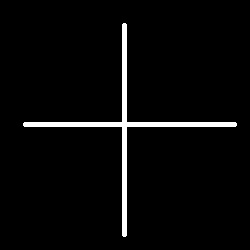}}~
    \subfloat[$f_1$]{\includegraphics[width=25mm]{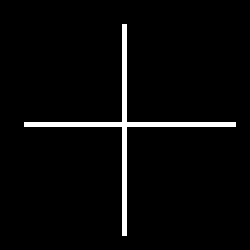}}~
    \subfloat[$f_2$]{\includegraphics[width=25mm]{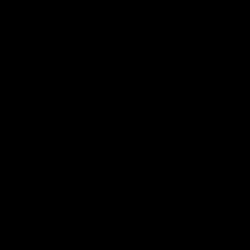}}~
    \subfloat[$f_3$]{\includegraphics[width=25mm]{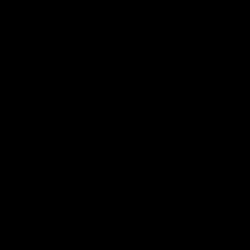}}~
    \subfloat[$f_4$]{\includegraphics[width=25mm]{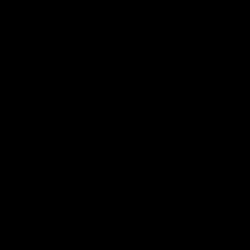}}\\
    \subfloat[$g = g_0$]{\includegraphics[width=25mm]{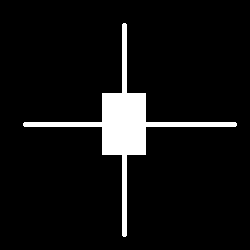}}~
    \subfloat[$g_1$]{\includegraphics[width=25mm]{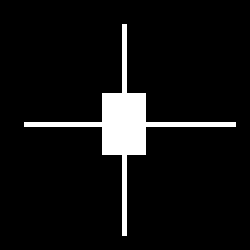}}~
    \subfloat[$g_2$]{\includegraphics[width=25mm]{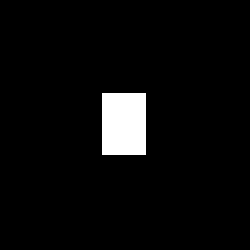}}~
    \subfloat[$g_3$]{\includegraphics[width=25mm]{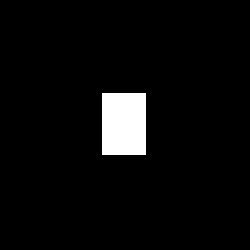}}~
    \subfloat[$g_4$]{\includegraphics[width=25mm]{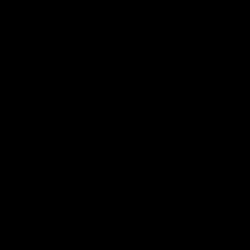}}
    \end{center}
    \caption{Two filtrations of 2-dimensional black pixels; that is, $f_0^{-1}(0) \subseteq f_1^{-1}(0) \subseteq f_2^{-1}(0) \subseteq f_3^{-1}(0) \subseteq f_4^{-1}(0)$ and $g_0^{-1}(0) \subseteq g_1^{-1}(0) \subseteq g_2^{-1}(0) \subseteq g_3^{-1}(0) \subseteq g_4^{-1}(0)$. Although images $f$ and $g$ share the same $1$-dimensional homology space $\bbZ_2$, the persistent homologies of these two images depict different lifespans. Indeed, the $1$-dimensional hole in (a)-(e) has the barcode $(0,2)$ while the hole in (f)-(j) has $(0,4)$.}
    \label{Fig. Opening filtration}
 \end{figure*}

Homology detects the hole structure in a given topological space, while it may omit some geometry of the based space. For example, two geometric objects with a single 1-dimensional hole in different sizes share the same first homology group (Figure \ref{Fig. Opening filtration}). As a generalization of homology, \textit{persistent homology} (PH) concerns sequences of topological spaces and their homologies. It was motivated by the works related to the Morse theory of Patrizio Frosini~\cite{frosini1992measuring} and Vanessa Robins~\cite{robins1999towards} in the 1990s. In Morse theory, a height function $f: M \rightarrow \bbR$ on a smooth manifold $M$ can form a sublevel set filtration of subspaces of $M$~\cite{EdelsbrunnerHarerbook2010}. The topological changes of such sublevel sets (e.g., the changes of Betti numbers) track the shape of $M$ along the direction of the height function and hence a descriptor (or fingerprint) of $M$. Persistent homology of height functions is now a well-known and fundamental tool in Morse theory and has many applications in theory~\cite{mischaikow2013morse,bubenik2010statistical,usher2016persistent} and data science~\cite{chung2018topological,delgado2014morse,gunther2011memory,kannan2019persistent}.

More generally, besides the smooth structures, suppose we have a sequence $X_1 \xrightarrow{f_1} X_2 \xrightarrow{f_2} \cdots \xrightarrow{f_{n-1}} X_n$ of topological spaces and continuous maps, then the functoriality of singular homology shown in Proposition \ref{Proposition: functoriality of singular homology} induces a sequence of homologies as follows:
\begin{equation*}
    H_{q}(X_1;R) \xrightarrow{H_{q}(f_1)} H_{q}(X_2;R) \xrightarrow{H_{q}(f_2)} \cdots \xrightarrow{H_{q}(f_{n-1})} H_{q}(X_n;R), 
\end{equation*}
where $q$ is an arbitrary non-negative integer, and $H_q(X_i)$, $H_q(f_i)$ are vector spaces and linear transformations over $\bbZ_2$. Because continuous maps can deform the geometry of spaces (e.g., sizes, lengths, and connectivity), the changes in homological cycles and Betti numbers depict how the hole structures in the spaces changed among the continuous deformation.

Computing homologies connected by continuous maps is challenging in real applications, so one usually considers a chain of filtered topological spaces with subspace relations. A tower of such topological spaces is called a \textit{filtration}. We list the formal definition of filtration as follows.
\begin{def.}[\cite{EdelsbrunnerHarerbook2010}]
A \textbf{filtration} of topological spaces is a sequence $\emptyset = X_0, X_1, X_2, ..., X_n$ of topological spaces such that $X_i$ is a subspace of $X_{i+1}$ for each $i \in \{ 0,1,...,n-1 \}$. We usually use the chain
\begin{equation*}
\mathcal{F}: \emptyset = X_0 \subseteq X_1 \subseteq X_2 \subseteq \cdots \subseteq X_n  
\end{equation*}
of topological spaces to denote a filtration of topological spaces. 
\end{def.}
Because $H_q(\cdot; R) : \mathfrak{Top} \rightarrow \mathfrak{Mod}_R$ is a functor, a filtration of topological spaces $\emptyset = X_0 \subseteq X_1 \subseteq \cdots \subseteq X_n$ and a non-negative integer $q \geq 0$ induce a sequence of $R$-modules and $R$-module homomorphisms:
\begin{equation}
\label{eq. Persistent Homology}
 {\rm PH}_q : 0 = H_q(\emptyset;R) \xrightarrow{\rho_{0,1}} H_q(X_1;R) \xrightarrow{\rho_{1,2}} H_q(X_2;R) \rightarrow \cdots \rightarrow H_q(X_n;R)
\end{equation}
where the $R$-module homomorphism $\rho_{i,j} : H_q(X_i;R) \rightarrow H_q(X_j;R)$ for $i \leq j$ is induced by the inclusion $X_i \hookrightarrow X_j$. Based on the functoriality of singular homology on the sequence \eqref{eq. Persistent Homology}, we define $\rho_{i,j} = \rho_{j - 1,j} \circ \rho_{j - 2,j - 1} \circ \cdots \circ \rho_{i, i+1}$ for every $i \leq j$ in $\{ 0,1, ..., n\}$, then the $\rho_{i,j}$ is also the $R$-module homomorphism induced by the inclusion map $X_i \hookrightarrow X_j$.
\begin{def.}
[\cite{EdelsbrunnerHarerbook2010}]
Suppose $\mathcal{F}: \emptyset = X_0 \subseteq X_1 \subseteq \cdots \subseteq X_n$ is a filtration of topological spaces. Then, for every ring $R$ and $q \in \bbZ_{\geq 0}$, we call the sequence defined in \eqref{eq. Persistent Homology} is the $\textbf{q}$\textbf{-}${\textbf{th}}$ \textbf{persistent homology} of the filtration $\mathcal{F}$.
\end{def.}
One of the primary purposes of persistent homology is to track the lifespans of local holes, i.e., the births/deaths of connected components, loops, and higher dimensional voids. To tackle this problem, H. Edelsbrunner and J. Harer proposed the \textit{persistence barcode} of persistent homology to detect such topological changes \cite{edelsbrunner2008persistent,EdelsbrunnerHarerbook2010}. We refer to the definition of persistence barcodes as follows.
\begin{def.}
[\cite{edelsbrunner2008persistent, EdelsbrunnerHarerbook2010}]
\label{Def. Persistence Barcodes}
Suppose $\emptyset = X_0 \subseteq X_1 \subseteq \cdots \subseteq X_n$ is a filtration of topological spaces and $0 \rightarrow H_q(X_1;F) \rightarrow \cdots \rightarrow H_q(X_n;F)$ is the induced $q^{\rm th}$ persistent homology over a field $F$. Let $s_i$ be an element in $H_q(X_i;F)$ ($i \geq 1$). Then we have the following definitions:
\begin{itemize}
    \item[\rm (a)] $s_i$ is said to \textbf{be born} at $i$ if $s_i \notin {\rm im}(\rho_{i-1,i})$, $i$ is called the \textbf{birth} of $s_i$;
    \item[\rm (b)] $s_i$ is said to \textbf{die} at $j$ if $\rho_{i,j-1}(s_i) \notin {\rm im}(\rho_{i-1,j-1})$ and $\rho_{i,j}(s_i) \in {\rm im}(\rho_{i-1,j})$, $j$ is called the \textbf{death} of $s_i$.
\end{itemize}
If $s_i$ is still alive at $n$, we define the death of $s_i$ to be $+\infty$ (up to this filtration). The tuple $(i,j)$ of $s_i$ is called the \textbf{persistence barcode} of the element $s_i \in H_q(X_i;F)$. The multiset of all persistence barcodes of non-repeated representative generators in all $H_q(X_i; F)$ is called the \textbf{persistence diagram} of the filtration. 
\end{def.}
For example, by considering the geometry of 2D black objects, rows in Figure \ref{Fig. Opening filtration} define two filtrations of subspaces in $\bbR^2$, and the induced  first persistent homologies (over $\bbZ_2$) are 
\begin{equation}
\label{Eq. Cross}
\begin{split}
\xymatrix@+0.0em{
				& \bbZ_2
				\ar[r]^{\id_{\bbZ_2}} 
                & \bbZ_2
                \ar[r] 
                & 0
				\ar[r] 
                & 0
                \ar[r] 
                & 0
                & \text{and}
                & \bbZ_2
				\ar[r]^{\id_{\bbZ_2}} 
                & \bbZ_2
                \ar[r]^{\id_{\bbZ_2}} 
                & \bbZ_2
				\ar[r]^{\id_{\bbZ_2}} 
                & \bbZ_2
                \ar[r] 
                & 0
}.
\end{split}
\end{equation}
By definition, the $1$-dimensional hole in Figure \ref{Fig. Opening filtration}(a)-(e) has the barcode $(0,2)$. On the other hand, the hole in Figure \ref{Fig. Opening filtration}(f)-(j) has barcode $(0,4)$.

There are many different ways to construct filtrations and compute their persistent homology. A typical one is the Vietoris--Rips complexes for the point-cloud data. For a (finite) set $\mathcal{X}$ in the $n$-dimensional Euclidean space $\bbR^n$ and a fixed positive real number $\epsilon > 0$, one explores the intersections of $n$-dimensional balls centered at points $x$ in $\mathcal{X}$ with radius $\epsilon$. Regarding points in $\mathcal{X}$ as the vertices of a simplicial complex, higher repeated regions lead to higher dimensional simplices in $\bbR^n$.  The strategy of the Vietoris--Rips complex is to enlarge the radius to construct a filtration of simplicial complexes~\cite{latschev2001vietoris,ghrist2005coverage,de2006coordinate}.

\begin{figure}
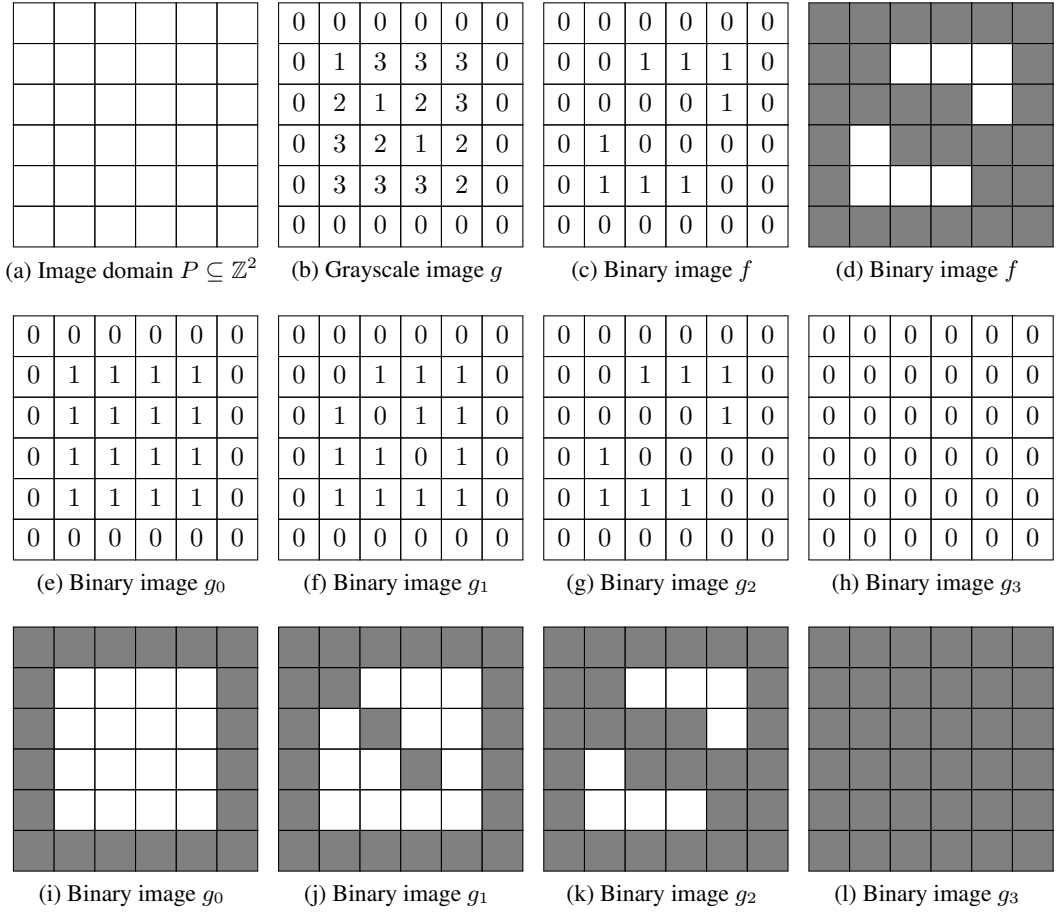
%
  \centering
  \subfloat[Image domain $P \subseteq \bbZ^2$]{
\begin{ytableau}
    \ &  \ & \ & \ & \ &\\
    & & & & & \\
    & & & & & \\
    & & & & & \\
    & & & & & \\
    & & & & & \\
\end{ytableau}}~
  \subfloat[Grayscale image $g$]{
\begin{ytableau}
    0 & 0 & 0 & 0 & 0 & 0 \\
    0 & 1 & 3 & 3 & 3 & 0\\
    0 & 2 & 1 & 2 & 3 & 0\\
    0 & 3 & 2 & 1 & 2 & 0\\
    0 & 3 & 3 & 3 & 2 & 0\\
    0 & 0 & 0 & 0 & 0 & 0 \\
\end{ytableau}}~
  \subfloat[Binary image $f$]{
  \begin{ytableau}
    0 & 0 & 0 & 0 & 0 & 0 \\
    0 & 0 & 1 & 1 & 1 & 0\\
    0 & 0 & 0 & 0 & 1 & 0\\
    0 & 1 & 0 & 0 & 0 & 0\\
    0 & 1 & 1 & 1 & 0 & 0\\
    0 & 0 & 0 & 0 & 0 & 0 \\
\end{ytableau}}~ 
  \subfloat[Binary image $f$]{
  \begin{ytableau}
   *(gray) & *(gray) & *(gray) & *(gray) & *(gray) & *(gray)\\
   *(gray) & *(gray) & & & & *(gray) \\
   *(gray) & *(gray) & *(gray) & *(gray) & & *(gray) \\
   *(gray) & & *(gray) & *(gray) & *(gray) & *(gray) \\
   *(gray) & & & & *(gray) & *(gray) \\
   *(gray) & *(gray) & *(gray) & *(gray) & *(gray) & *(gray)\\
  \end{ytableau}}\\%
  \subfloat[Binary image $g_0$]{
  \begin{ytableau}
    0 & 0 & 0 & 0 & 0 & 0 \\
    0 & 1 & 1 & 1 & 1 & 0\\
    0 & 1 & 1 & 1 & 1 & 0\\
    0 & 1 & 1 & 1 & 1 & 0\\
    0 & 1 & 1 & 1 & 1 & 0\\
    0 & 0 & 0 & 0 & 0 & 0 \\
\end{ytableau}}~
  \subfloat[Binary image $g_1$]{
  \begin{ytableau}
    0 & 0 & 0 & 0 & 0 & 0 \\
    0 & 0 & 1 & 1 & 1 & 0\\
    0 & 1 & 0 & 1 & 1 & 0\\
    0 & 1 & 1 & 0 & 1 & 0\\
    0 & 1 & 1 & 1 & 1 & 0\\
    0 & 0 & 0 & 0 & 0 & 0 \\
\end{ytableau}}~
\subfloat[Binary image $g_2$]{
  \begin{ytableau}
    0 & 0 & 0 & 0 & 0 & 0 \\
    0 & 0 & 1 & 1 & 1 & 0\\
    0 & 0 & 0 & 0 & 1 & 0\\
    0 & 1 & 0 & 0 & 0 & 0\\
    0 & 1 & 1 & 1 & 0 & 0\\
    0 & 0 & 0 & 0 & 0 & 0 \\
\end{ytableau}}~
\subfloat[Binary image $g_3$]{
  \begin{ytableau}
    0 & 0 & 0 & 0 & 0 & 0 \\
    0 & 0 & 0 & 0 & 0 & 0\\
    0 & 0 & 0 & 0 & 0 & 0\\
    0 & 0 & 0 & 0 & 0 & 0\\
    0 & 0 & 0 & 0 & 0 & 0\\
    0 & 0 & 0 & 0 & 0 & 0 \\
\end{ytableau}}\\
\subfloat[Binary image $g_0$]{
  \begin{ytableau}
   *(gray) & *(gray) & *(gray) & *(gray) & *(gray) & *(gray)\\
   *(gray) & & & & & *(gray) \\
   *(gray) & & & & & *(gray) \\
   *(gray) & & & & & *(gray) \\
   *(gray) & & & & & *(gray) \\
   *(gray) & *(gray) & *(gray) & *(gray) & *(gray) & *(gray)\\
  \end{ytableau}}~
\subfloat[Binary image $g_1$]{
  \begin{ytableau}
   *(gray) & *(gray) & *(gray) & *(gray) & *(gray) & *(gray)\\
   *(gray) & *(gray) & & & & *(gray) \\
   *(gray) &  & *(gray) & & & *(gray) \\
   *(gray) & & & *(gray) & & *(gray) \\
   *(gray) & & & & & *(gray) \\
   *(gray) & *(gray) & *(gray) & *(gray) & *(gray) & *(gray)\\
  \end{ytableau}}~
\subfloat[Binary image $g_2$]{
  \begin{ytableau}
   *(gray) & *(gray) & *(gray) & *(gray) & *(gray) & *(gray)\\
   *(gray) & *(gray) & & & & *(gray) \\
   *(gray) & *(gray) & *(gray) & *(gray) & & *(gray) \\
   *(gray) & & *(gray) & *(gray) & *(gray) & *(gray) \\
   *(gray) & & & & *(gray) & *(gray) \\
   *(gray) & *(gray) & *(gray) & *(gray) & *(gray) & *(gray)\\
  \end{ytableau}}~
\subfloat[Binary image $g_3$]{
  \begin{ytableau}
   *(gray) & *(gray) & *(gray) & *(gray) & *(gray) & *(gray)\\
   *(gray) & *(gray) & *(gray) & *(gray) & *(gray) & *(gray)\\
   *(gray) & *(gray) & *(gray) & *(gray) & *(gray) & *(gray)\\
   *(gray) & *(gray) & *(gray) & *(gray) & *(gray) & *(gray)\\
   *(gray) & *(gray) & *(gray) & *(gray) & *(gray) & *(gray)\\
   *(gray) & *(gray) & *(gray) & *(gray) & *(gray) & *(gray)\\
  \end{ytableau}}
  \caption{First row: a $6 \times 6$ image domain $P$ in $\bbZ^2$, a grayscale image $g: P \rightarrow \{ 0,1,2,3\}$, and a binary image $f: P \rightarrow \{ 0,1 \}$. Figures (c) and (d) are two different representations for the image $f$. In a binary image $f$ as in (d), pixels with a value of $0$ represent the black pixels of the image. Second row: a filtration of binary images made by image $g$ and thresholds $0, 1, 2,$ and $3$. Third row: the white-black pixel representations of images in the second row.}
  \label{Fig. Definition of images}
\end{figure}

As shown in Figure \ref{Fig. Opening filtration} and equation~\eqref{Eq. Cross}, except for the point-cloud data,  one can also construct filtrations of digital images and compute their persistent homology. We referred to an $m$-dimensional digital image as a function $f: P \xrightarrow{} \bbR_{\geq 0}$ from a non-empty set $P$ of $\bbZ^m$ to the set of all non-negative real numbers (cf.~\cite{chung2022multi}). An image $f$ is called \textit{binary} if its range is contained in the binary set $\{ 0,1 \}$ and called \textit{grayscale} for otherwise. For a binary image, the primage of zero $f^{-1}(0)$ referred to the set of all \textit{black pixels} of $f$, and $f^{-1}(1)$ denotes the set of all \textit{white pixels} of $f$. Viewing each black pixel as a closed cube in $\bbR^m$, we regard $f^{-1}(0)$ as a subspace of $\bbR^m$ and consider its topological properties. The first row in Figure \ref{Fig. Definition of images} provides examples of $2$-dimensional grayscale and binary digital images.

As in Figure \ref{Fig. Opening filtration}, one can construct filtrations of images by using image processing techniques on a given binary one. Another typical method of building filtrations is operating the sub-level sets of a grayscale image. For a image $f: P \xrightarrow{} \bbR_{\geq 0}$ and a threshold $t \in \bbR$, we define a binary image $f_t : P \xrightarrow{} \{ 0,1 \}$ by setting $f_t(x) = 0$ if $f(x) \leq t$ and $f_t(x) = 1$ for otherwise. Then $f_{t_1}^{-1}(0) \subseteq f_{t_2}^{-1}(0) \subseteq \cdots \subseteq f_{t_n}^{-1}(0)$ for $t_1 \leq t_2 \leq \cdots \leq t_n$. The second the third row in Figure \ref{Fig. Definition of images} illustrate how sub-level sets of a grayscale image form a filtration of black pixels. In particular, the $0$-th and $1$-th persistence diagrams of the filtrations are $\{ (0,+\infty) \}$ and $\{ (0,3), (2,3) \}$. For readers who are interested in persistent homology on digital images, see~\cite{kaczynski2004computational} for more information.

In this paper, we focus on 2-dimensional binary images and their local homology. We combine image segmentation techniques, local homology, and persistence barcodes to illustrate how to estimate and detect the positions of holes in 2D binary images. The combination of this detection method with more image processing techniques (such as mathematical morphology and sub-level set filtration) will be our future work.

\section{Our Approaches}
\label{Sec. Our Approaches}

The section is separated into three parts. First, we quote the definitions of local systems and short persistent homology in our previous work~\cite{hu2021sheaf}. \textit{Local systems} and \textit{short persistent homology} induce a cellular sheaf structure of topological spaces and can depict the spatial merging relations via their global/local sections~\cite{hu2021sheaf}. We discuss the relationship between hole positions, global/local sections, and persistence barcodes on local systems  (Section \ref{Subsec. Persistent Homology of Local Systems}). Second, we introduce how we adapt the theory to digital images and implement the method (Section \ref{Subsec. Local Systems in Binary Images}). Finally, we discuss some properties of the proposed framework, such as the relationship between local systems, the location of holes, and image noises (Section \ref{Subsec. Discussion}).

\subsection{Persistent Homology of Local Systems}
\label{Subsec. Persistent Homology of Local Systems}

For a topological space $X$ and a concerned local region $A$ of $X$, the \textit{relative homology} $H_q(X, A)$ considers the equivalence classes of cycles in $X$ that do not meet the subspace $A$. One can formulate the relative homology of $X$ and $A$ by $H_q(X, A) = Z_q(X,A)/B_q(X,A)$, where $Z_q(X,A) = \{ c \in S_q(X) : \partial_q(c) \in S_{q-1}(A) \} = \partial_q^{-1}(S_q(A))$ is the set of all chains in $S_{q}(X)$ with boundaries in $S_{q-1}(A)$, and $B_q(X,A) = B_q(X) + S_q(A)$ is the submodule generated by all $q$-boundaries of $X$ and $q$-chains in $A$. Elements in $Z_q(X, A)$ and $B_q(X, A)$ are called relative $q$-cycles and relative $q$-boundaries of $X$, respectively~\cite{Greenberg}. 


Roughly speaking, relative homology detects holes in $X$ except for holes that are totally contained in $A$. More precisely, one can apply the snake lemma on the short exact sequence $0 \xrightarrow{} S_{\bullet}(A) \xrightarrow{\iota_\bullet} S_{\bullet}(X) \xrightarrow{\pi_\bullet} S_{\bullet}(X)/S_{\bullet}(A) \xrightarrow{} 0$ with canonical inclusion and projection to obtain the long exact sequence
\begin{equation}
\label{Eq. Long exact sequence of relative homology}
    \cdots \xrightarrow{ } H_q(A) \xrightarrow{\overline{\iota}_q} H_q(X) \xrightarrow{\overline{\pi}_q} H_q(X,A) \xrightarrow{\delta_q} H_{q-1}(A) \xrightarrow{\overline{\iota}_{q-1}} H_{q-1}(X) \xrightarrow{\overline{\pi}_{q-1}} H_{q-1}(X,A) \xrightarrow{ } \cdots.
\end{equation}
One can use barcode representation to detect hole structures in the spaces $H_\bullet(A), H_\bullet(X),$ and $ H_\bullet(X, A)$. For example, a non-zero element in $H_q(X) \setminus \im(\overline{\iota}_q)$ represents a hole in $X$ that is not totally emerged in the region $A$. On the other hand,  $c \in H_q(X)$ dies at $H_q(X,A)$ if the cycle $c$ does not represent a hole in $A$.
\begin{remark}
Because the long exact sequence in ~\eqref{Eq. Long exact sequence of relative homology} is a chain complex, every lifespan $b - d$ of a barcode $(b,d)$ is $1$. 
\end{remark}

Relative homology can capture holes contributed by $A$, $X \setminus A$, or both. However, it is difficult and expensive to implement and compute due to the complicated data representation. This paper proposes a relatively efficient method to detect hole relations and positions via persistent homology. To achieve this goal, we introduce here two main ideas proposed in our previous work, called \textit{local system} and \textit{short filtration}~\cite{hu2021sheaf}.

\begin{def.}[\cite{hu2021sheaf}]
Let $X$ be a topological space and $X_1, X_2$ be subspaces of $X$. The triad $(X, X_1, X_2)$ is called a \textbf{local system} (or an \textbf{admissible triad}) if $\cl_X(X_1) \cap \cl_X(X_2) = \emptyset$.   
\end{def.}

\begin{figure*}
    \begin{center}
    \includegraphics[width=150mm]{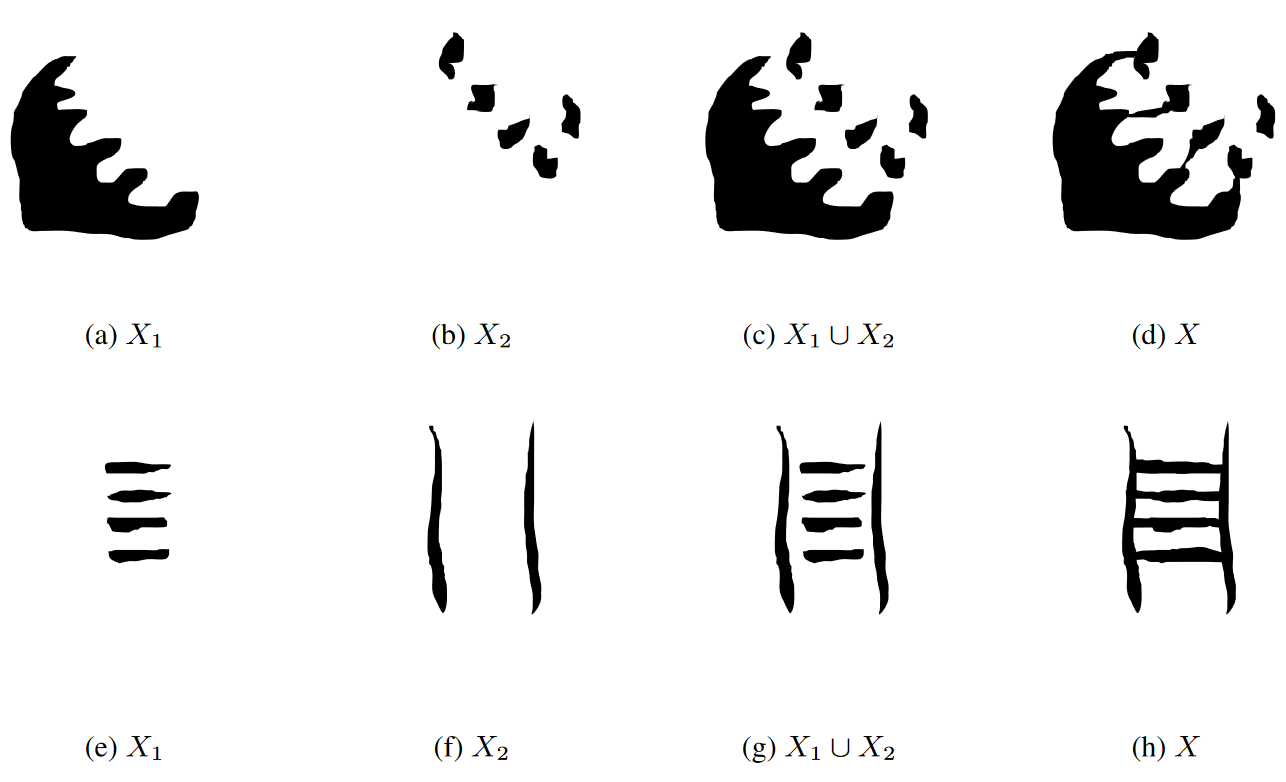}
    \end{center}
    \caption{Two local systems $(X,X_1,X_2)$ of topological spaces. Let $\Gamma_0$ denote the global section space of the sheaf structure $H_0(X_1) \xrightarrow{} H_0(X) \xleftarrow{} H_0(X_2)$. Then we have the following information. The first row: $H_0(X_1) \simeq \bbZ_2$, $H_0(X_2) \simeq \bbZ_2^5$, $H_0(X_1 \cup X_2) \simeq \bbZ_2^6$, $H_0(X) \simeq \bbZ_2^2$, and $\Gamma_0 \simeq \bbZ_2^4$. The second row: $H_0(X_1) \simeq \bbZ_2^4$, $H_0(X_2) \simeq \bbZ_2^2$, $H_0(X_1 \cup X_2) \simeq \bbZ_2^6$, $H_0(X) \simeq \bbZ_2$, and $\Gamma_0 \simeq \bbZ_2^5$.}
    \label{Fig. Sheaf Structures}
 \end{figure*}

For any topological space $X$ and its subspaces $X_1$ and $X_2$, we have the following definition.

\begin{def.}[\cite{hu2021sheaf}]
\label{Def. Short filtration}
Let $(X, X_1, X_2)$ be a triad of topological spaces with $X_1 \subseteq X$ and $X_2 \subseteq X$. This triad leads to two filtrations $\emptyset \subseteq X_1 \subseteq X_1 \cup X_2 \subseteq X$ and $\emptyset \subseteq X_2 \subseteq X_1 \cup X_2 \subseteq X$. We call them \textbf{short filtrations} of the triad $(X, X_1, X_2)$.
\end{def.}

Focusing on the first one in Definition \ref{Def. Short filtration}, the birth information at $H_\bullet(X_1 \cup X_2; F)$ depicts whether $X_2$ contains a homological generator that cannot be represented via generators in $H_\bullet(X_1; F)$.  When $\cl_X(X_1) \cap \cl_X(X_2) = \emptyset$, the homology $H_\bullet(X_1 \cup X_2; F)$ is canonically isomorphic to the space $H_\bullet(X_1; F) \oplus H_\bullet(X_2; F)$ since $X_1$ and $X_2$ are two path-connected components of $X_1 \cup X_2$. In this case, every generator $s_2$ in $H_\bullet(X_2; F)$ is born at $H_\bullet(X_1 \cup X_2; F)$ of the persistent homology $0 \xrightarrow{} H_\bullet(X_1; F) \xrightarrow{} H_\bullet(X_1 \cup X_2; F) \xrightarrow{} H_\bullet(X; F)$ and dies at $H_\bullet(X; F)$ if there is an $s_1 \in H_\bullet(X_1; F)$ such that $s_1$ and $s_2$ represent the same homological generator in $H_\bullet(X; F)$. This property will benefit computing the homological changes of holes in $X_1$, $X_2$, and $X$. Furthermore, we will show in Section \ref{Subsec. Local Systems in Binary Images} that the condition $\cl_X(X_1) \cap \cl_X(X_2) = \emptyset$ can be easily established in image data through elementary image processing techniques.

In \cite{hu2021sheaf}, we applied the two filtrations of a local system $(X, X_1, X_2)$ to construct the following \textit{cellular sheaf} structure: 
\begin{equation*}
\xymatrix@+1.0em{
                & H_{q}(X_1;F)
				\ar[r]^{\rho_1}
                & H_{q}(X;F)
                \\
                & 
        	& H_{q}(X_2;F)
                \ar[u]_{\rho_2}
}    
\end{equation*}
where $q$ is any non-negative integer, $F$ is a fixed field, and $\rho_1, \rho_2$ are the $F$-linear transformations induced by the inclusions $X_1 \hookrightarrow X$ and $X_2 \hookrightarrow X$. We often call the maps $\rho_1, \rho_2$ \textit{restriction maps}. A pair $(s_1,s_2) \in H_{q}(X_1;F) \oplus H_{q}(X_2;F)$ is called a \textit{global} \textit{section} of the sheaf if $\rho_1(s_1) = \rho_2(s_2)$. We use $\Gamma$ to denote the subspace of all global sections in $H_{q}(X_1;F) \oplus H_{q}(X_2;F)$, and it can be sculptured by the following theorem.

\begin{theorem}
\label{Theorem: Global section space of the simplest sheaf}
For the following sheaf structure of $F$-vector spaces and $F$-linear maps:
\begin{equation*}
\xymatrix@-1.5em{
V \ar[rr]^{f} & & P & \\
 & & & , \\
 & & W \ar[uu]_{g}& \\
}
\end{equation*}
we define $\phi : V \oplus W \rightarrow P$ by $(v,w) \longmapsto f(v) - g(w)$. Then, $\phi$ is also an $F$-linear linear map and $(V \oplus W) / \Gamma \simeq {\rm im}(\phi)$, where $\Gamma = \{ (v,w) : f(v) = g(w) \}$ is the space of global sections.

In particular, $\dim(\Gamma) = \dim(V) + \dim(W) - \dim(\im(\phi))$ if the spaces $V, W, P$ are finite-dimensional. In addition, $\dim(\Gamma) = \dim(V) + \dim(W) - \dim(P)$ if $\phi$ is onto.
\end{theorem}
\begin{proof}
It is evident that $\phi$ is $F$-linear. Because $(v,w) \in \ker(\phi)$ if and only if $f(v) = f(w)$. By the first isomorphism theorem of modules, the theorem follows.    
\end{proof}

Let $(X, X_1, X_2)$ be a local system of topological spaces and $\Gamma$ the global section space of the sheaf structure $H_q(X_1;F) \xrightarrow{} H_q(X;F) \xleftarrow{} H_q(X_2;F)$. Examples shown in Figure \ref{Fig. Sheaf Structures} depict that the vector spaces $H_q(X_1;F)$, $H_q(X_2;F)$, $H_q(X_1 \cup X_2;F)$, $H_q(X;F)$, and $\Gamma$ can be totally different. In other words, the global section space provides an additional than the homology of $X_1, X_2, X_1 \cup X_2,$ and $X$. Actually, suppose we have a sequence $(X_i, X_{i1}, X_{i2})$ of local systems that satisfy $X_{i1} \subseteq X_{(i+1)1}$, $X_{i2} \subseteq X_{(i+1)2}$, and $X_{i} \subseteq X_{i+1}$, then we have the following commutative diagram:
\begin{equation*}
    \xymatrix@+0.0em{
H_q(X_{11};F) \ar[d]_{\rm res.} \ar[r]_{\rm \phi_{11}} & H_q(X_{21};F) \ar[d]_{\rm res.} \ar[r]_{\rm \phi_{21}} & H_q(X_{31};F) \ar[r] \ar[d]_{\rm res.} & \cdots\\
H_q(X_1;F)    \ar[r]_{\rm \phi_{1}}                    & H_q(X_2;F) \ar[r]_{\rm \phi_{2}} & H_q(X_{3};F) \ar[r] & \cdots\\
H_q(X_{12};F) \ar[u]^{\rm res.} \ar[r]_{\rm \phi_{12}} & H_q(X_{22};F) \ar[u]^{\rm res.} \ar[r]_{\rm \phi_{22}} & H_q(X_{32};F) \ar[r] \ar[u]^{\rm res.} & \cdots
}
\end{equation*}
where $\phi_{ij}$ and $\phi_i$ are the $F$-linear maps induced by the inclusions. One can check  the following sequence is also valid:
\begin{equation*}
    \xymatrix@+3.5em{
\Gamma_1 \ar[r]^{\phi_{11} \oplus \phi_{12}|_{\Gamma_1}}                    & \Gamma_2 \ar[r]^{\phi_{21} \oplus \phi_{22}|_{\Gamma_2}} & \Gamma_3 \ar[r] & \cdots\\
}.
\end{equation*}
In other words, except for computing single global section spaces, one can also consider the persistent homology of global section spaces induced by any filtered local systems of topological spaces.

Theorem \ref{Theorem: Global section space of the simplest sheaf} presents a way to compute global section spaces. However, on many occasions, computing the image of $\phi$ in the theorem may be infeasible. To tackle this, we previously  proposed an approximation method using persistent
homology~\cite{hu2021sheaf}. We quote the method as the following theorem.

\begin{theorem}[Theorem 2.3.1~\cite{hu2021sheaf}]
\label{Theorem: Local section approximation}
Let $R$ be a commutative ring with identity. Let $(X, X_1, X_2)$ be a local system of topological spaces and $q$ a non-negative integer. Let $\mathcal{G}_1$ be the short filtration $\emptyset \subseteq X_1 \subseteq X_1 \cup X_2 \subseteq X$ and $s_2 \in H_q(X_2;R)$ a non-zero element. Then the followings are equivalent:
\begin{itemize}
    \item[\rm (a)] There is an $s_1 \in H_q(X_1;R)$ such that $(s_1, s_2) \in H_q(X_1;R) \oplus H_q(X_2;R)$ is global section;
    \item[\rm (b)] $\widetilde{s_2} := \omega_2(s_2)$ has barcode $(2,3)$ in the PH $\mathcal{P}_q(\mathcal{G}_1) : 0 \xrightarrow{} H_q(X_1;R) \xrightarrow{} H_q(X_1 \cup X_2;R) \xrightarrow{} H_q(X;R)$.
\end{itemize}
\end{theorem}

For a local system $(X,X_1,X_2)$, numbers of barcode $(2,3)$ in $\mathcal{P}_q(\mathcal{G}_1)$ records how many homological non-zero generators in $H_q(X_2;R)$ that merge to a generator in $H_q(X_1;R)$. In~\cite{hu2021sheaf}, we defined it as the \textit{$q$-th local merging number}.

\begin{def.}[\cite{hu2021sheaf}]
Let $(X, X_1, X_2)$ be a local system of topological spaces and $q \geq 0$. We define the $\textbf{q}$\textbf{-th local merging number} of $X_1$ and $X_2$ as the numbers of barcodes $(2,3)$ in $\mathcal{P}_q(\mathcal{G}_1)$ and denote it by $m_q(X_1;X_2)$.   
\end{def.}

We use the two pairs in Figure \ref{Fig. Sheaf Structures} to explain the local merging numbers. For the first row, we have $m_0(X_1; X_2) = 5$ since there are $5$ connected components that merge to $X_1$. On the other hand, $m_0(X_2; X_1) = 5$. Similarly, the $m_0(X_1; X_2)$ and $m_0(X_2; X_1)$ of the second row are $2$ and $4$, respectively. In particular, these two examples show the local merging numbers $m_0(X_1; X_2)$ and $m_0(X_2; X_1)$ are not equal in general. Actually, one can prove that $\max\{ m_0(X_1;X_2), m_0(X_2;X_1) \} \leq \dim(\Gamma) \leq m_0(X_1;X_2) + m_0(X_2;X_1)$~\cite{HUphdthesis}.

When $q = 0$, the local merging numbers $m_0(X_1; X_2)$ records how many connected components in $X_2$ connect to components in $X_1$ synchronously. In our previous work, we show that local regions with high $0$-local merging numbers are likely to be more joint parts of the ambient space and have the potential to analyze handwritten text with texture data~\cite{hu2021sheaf, HUphdthesis}. In these works, we focus on local merging numbers in dimension $1$ and $(2,3)$ barcodes in short filtrations, while the geometric meanings of higher dimensional merging numbers and $(3,+\infty)$ barcodes are still unknown.

In the following theorem, we show that the number of barcodes $(3,+\infty)$ in a short filtration can verify whether $X_1$ and $X_2$ contribute a hole (with dimension $\geq 1$) in $X$.

\begin{theorem}
\label{Theroem for barcode (3,inf)}
Let $F$ be a field. Let $(X, X_1, X_2)$ be a local system of topological spaces and $q$ a non-negative integer.  Let $\Gamma$ be the global section space of the sheaf structure $H_q(X_1;F) \xrightarrow{} H_q(X;F) \xleftarrow{} H_q(X_2;F)$. Then the number of $(3,+\infty)$ in the PH $\mathcal{P}_q(\mathcal{G}_1) : 0 \xrightarrow{} H_q(X_1;F) \xrightarrow{} H_q(X_1 \cup X_2;F) \xrightarrow{} H_q(X;F)$ equals
\begin{equation*}
    \dim_F(H_q(X;F)) - \dim_F(H_q(X_1;F)) - \dim_F(H_q(X_2;F)) + \dim_F(\Gamma). 
\end{equation*}
\end{theorem}
\begin{proof}
Let $\rho_1: H_q(X_1; F) \xrightarrow{} H_q(X; F)$ and $\rho_2: H_q(X_2; F) \xrightarrow{} H_q(X; F)$ be the canonical linear transformations that are induced by the inclusions.  Define $\phi = \rho_1 - \rho_2: H_q(X_1) \oplus H_q(X_2) \xrightarrow{} H_q(X)$, then 
\begin{equation}
\label{Eq.-1: (3,inf) barcode}
\dim_F(\Gamma) = \dim_F(H_q(X_1; F)) + \dim_F(H_q(X_2; F)) - \dim_F(\im(\phi))    
\end{equation}
by Theorem \ref{Theorem: Global section space of the simplest sheaf}. Because $H_q(X_1 \cup X_2; F)$ is canonically isomorphic to $H_q(X_1; F) \oplus H_q(X_2; F)$, the images of $\rho_1 - \rho_2$ and the map $H_q(X_1 \cup X_2; F) \xrightarrow{} H_q(X; F)$ in $\mathcal{P}_q(\mathcal{G}_1)$ are equal. Then the number of barcodes $(3,+\infty)$ in the persistent homology $\mathcal{P}_q(\mathcal{G}_1)$ counts the dimension of the space $H_q(X;F)/\im(\phi)$. Therefore, 
\begin{equation}
\label{Eq.-2: (3,inf) barcode}
\begin{split}
\# \{ \text{barcode } (3,+\infty) \text{ in } \mathcal{P}_q(\mathcal{G}_1) \} &= \dim_F(H_q(X;F)) - \dim_F(\im(\phi)).    
\end{split}
\end{equation}
By plugging equation \eqref{Eq.-1: (3,inf) barcode} into equation \eqref{Eq.-2: (3,inf) barcode}, the theorem follows.
\end{proof}

If $c \in Z_q(X_1;F) \subseteq Z_q(X;F)$ is a $q$-cycle that represents a q-dimensional hole in $X_1$, then $c$ must have a barcode $(1,\star)$ in the persistent homology $\mathcal{P}_q(\mathcal{G}_1)$. On the other hand, $c \in Z_q(X_2;F) \subseteq Z_q(X;F)$ representing a hole in $X_2$ implies that it has a barcode $(2,\star)$. In other words, the number of barcodes $(3,+\infty)$ in the persistent homology $\mathcal{P}_q(\mathcal{G}_1)$
records how many $q$-holes in $X$ are “supported” by both $X_1$ and $X_2$. In particular, removing either $X_1$ or $X_2$ will make those holes disappear. Intuitively, those holes are constructed by gluing the parts by $X_1$ and $X_2$, and hence we have the following definition.

\begin{def.}
Let $(X, X_1, X_2)$ be a local system of topological spaces and $q \geq 0$. We define the $\textbf{q}$\textbf{-th local outer-merging number} of $X_1$ and $X_2$ as the numbers of barcodes $(3,+\infty)$ in $\mathcal{P}_q(\mathcal{G}_1)$ and denote it by $o_q(X_1;X_2)$.   
\end{def.}

From the above discussion, it can be seen that the local outer-merging number records the contribution of a specific local area in the topological space to the hole structure. We present local outer-merging numbers for digital images in the next section (Section \ref{Subsec. Local Systems in Binary Images}). In addition, we will analyze the location of holes in the image by segmenting the image and the local outer-merging number of the corresponding region.

\begin{figure}
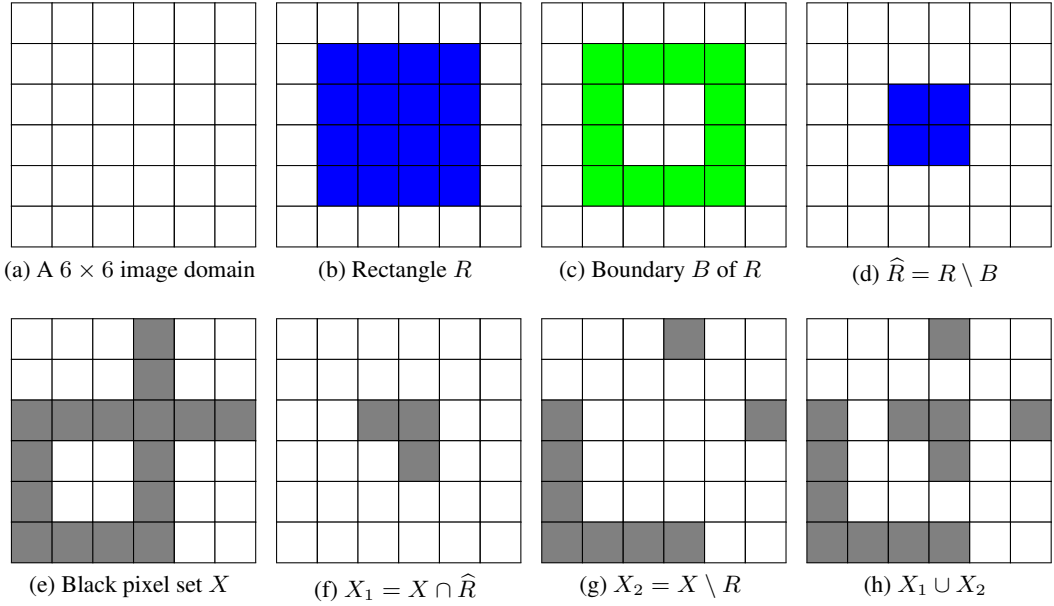
%
  \centering
  \subfloat[A $6 \times 6$ image domain]{
\begin{ytableau}
*(white) & *(white) & *(white) & *(white) & *(white) & *(white)\\
*(white) & *(white) & *(white) & *(white) & *(white) & *(white)\\
*(white) & *(white) & *(white) & *(white) & *(white) & *(white)\\
*(white) & *(white) & *(white) & *(white) & *(white) & *(white)\\
*(white) & *(white) & *(white) & *(white) & *(white) & *(white)\\
*(white) & *(white) & *(white) & *(white) & *(white) & *(white)\\
  \end{ytableau}}~
  \subfloat[Rectangle $R$]{
\begin{ytableau}
*(white) & *(white) & *(white) & *(white) & *(white) & *(white)\\
*(white) & *(blue) & *(blue) & *(blue) & *(blue) & *(white)\\
*(white) & *(blue) & *(blue) & *(blue) & *(blue) & *(white)\\
*(white) & *(blue) & *(blue) & *(blue) & *(blue) & *(white)\\
*(white) & *(blue) & *(blue) & *(blue) & *(blue) & *(white)\\
*(white) & *(white) & *(white) & *(white) & *(white) & *(white)\\
  \end{ytableau}}~
  \subfloat[Boundary $B$ of $R$]{
\begin{ytableau}
*(white) & *(white) & *(white) & *(white) & *(white) & *(white)\\
*(white) & *(green) & *(green) & *(green) & *(green) & *(white)\\
*(white) & *(green) & *(white) & *(white) & *(green) & *(white)\\
*(white) & *(green) & *(white) & *(white) & *(green) & *(white)\\
*(white) & *(green) & *(green) & *(green) & *(green) & *(white)\\
*(white) & *(white) & *(white) & *(white) & *(white) & *(white)\\
  \end{ytableau}}~
  \subfloat[$\widehat{R} = R \setminus B$]{
\begin{ytableau}
*(white) & *(white) & *(white) & *(white) & *(white) & *(white)\\
*(white) & *(white) & *(white) & *(white) & *(white) & *(white)\\
*(white) & *(white) & *(blue) & *(blue) & *(white) & *(white)\\
*(white) & *(white) & *(blue) & *(blue) & *(white) & *(white)\\
*(white) & *(white) & *(white) & *(white) & *(white) & *(white)\\
*(white) & *(white) & *(white) & *(white) & *(white) & *(white)\\
  \end{ytableau}} \\
  \subfloat[Black pixel set $X$]{
  \begin{ytableau}
*(white) & *(white) & *(white) & *(gray) & *(white) & *(white)\\
*(white) & *(white) & *(white) & *(gray) & *(white) & *(white)\\
*(gray) & *(gray) & *(gray) & *(gray) & *(gray) & *(gray)\\
*(gray) & *(white) & *(white) & *(gray) & *(white) & *(white)\\
*(gray) & *(white) & *(white) & *(gray) & *(white) & *(white)\\
*(gray) & *(gray) & *(gray) & *(gray) & *(white) & *(white)\\
  \end{ytableau}}~
  \subfloat[$X_1 = X \cap \widehat{R}$]{
\begin{ytableau}
*(white) & *(white) & *(white) & *(white) & *(white) & *(white)\\
*(white) & *(white) & *(white) & *(white) & *(white) & *(white)\\
*(white) & *(white) & *(gray) & *(gray) & *(white) & *(white)\\
*(white) & *(white) & *(white) & *(gray) & *(white) & *(white)\\
*(white) & *(white) & *(white) & *(white) & *(white) & *(white)\\
*(white) & *(white) & *(white) & *(white) & *(white) & *(white)\\
\end{ytableau}}~
  \subfloat[$X_2 = X \setminus R$]{
\begin{ytableau}
*(white) & *(white) & *(white) & *(gray) & *(white) & *(white)\\
*(white) & *(white) & *(white) & *(white) & *(white) & *(white)\\
*(gray) & *(white) & *(white) & *(white) & *(white) & *(gray)\\
*(gray) & *(white) & *(white) & *(white) & *(white) & *(white)\\
*(gray) & *(white) & *(white) & *(white) & *(white) & *(white)\\
*(gray) & *(gray) & *(gray) & *(gray) & *(white) & *(white)\\
  \end{ytableau}}~
\subfloat[$X_1 \cup X_2$]{
  \begin{ytableau}
*(white) & *(white) & *(white) & *(gray) & *(white) & *(white)\\
*(white) & *(white) & *(white) & *(white) & *(white) & *(white)\\
*(gray) & *(white) & *(gray) & *(gray) & *(white) & *(gray)\\
*(gray) & *(white) & *(white) & *(gray) & *(white) & *(white)\\
*(gray) & *(white) & *(white) & *(white) & *(white) & *(white)\\
*(gray) & *(gray) & *(gray) & *(gray) & *(white) & *(white)\\
  \end{ytableau}}
  \caption{An illustration of the construction of a local system in a 2D binary image. In this example, we have $m_0(X_1;X_2) = 3$, $o_0(X_1;X_2) = 0$, $m_1(X_1;X_2) = 0$, and $o_1(X_1;X_2) = 1$.}
  \label{Fig. Local section in images}
\end{figure}

\subsection{Local Systems in Binary Images}
\label{Subsec. Local Systems in Binary Images}

Section \ref{Subsec. Persistent Homology of Local Systems} introduces the local system and its persistent homology. Theorem \ref{Theorem: Local section approximation} and Theorem \ref{Theroem for barcode (3,inf)} tell us that counting the numbers of barcodes $(2,3)$ and $(3,+\infty)$ in $(X, X_1, X_2)$ can detect the glue relationship of local objects in $X$. Among them, constructing the admissible triad $(X, X_1, X_2)$ is the most crucial part of the calculation. For an object $X$ in $\bbR^n$ and a bounded $A \subseteq X$, one can choose $r_1, r_2 > 0$ with $r_1 < r_2$ such that $A \subseteq \mathbf{B}(\mathbf{0},r_1)$ and define $X_2 = X \cap \{ \bfx \in \bbR^n : |\bfx| \geq  r_2 \}$. Then, $\cl_X(X_1) \cap \cl_X(X_2) = \emptyset$. Based on the same idea, the section presents a more efficient way to build local systems in binary images.

As we introduced in Section \ref{Subsec. Persistent Homology}, a $2$-dimensional image is identified as a non-negative real-valued function $f : P \xrightarrow{} \bbR_{\geq 0}$ on a discrete 2D rectangle $P = ([a,b] \times [c,d]) \cap \mathbb{Z}^2$, where $a, b, c, d$ are integers with $a \leq b$ and $c \leq d$. In the paper, we focus on the geometric realization of black pixels of a binary image and compute its homology (see Figure \ref{Fig. Definition of images}(d) and Figure \ref{Fig. Local section in images}). For a binary image $f: P \xrightarrow{} \{ 0,1 \}$, we consider the black pixel set $f^{-1}(0)$ and denote it by $X \subseteq P$. We use a rectangle in $P$ to cover a concerned region of $X$, say $R = ([a_1,b_1] \times [c_1,d_1]) \cap \mathbb{Z}^2$ with $a \leq a_1 \leq b_1 \leq b$ and $c \leq c_1 \leq d_1 \leq d$ (see Figure \ref{Fig. Local section in images}(b)). Consider \begin{equation*}
B = \bbZ^2 \cap \bigg( (\{ a_1 \} \times [c_1, d_1]) \cup (\{ b_1 \} \times [c_1, d_1]) \cup ([a_1, b_1] \times \{ c_1 \})  \cup ([a_1, b_1] \times \{ d_1 \})  \bigg)   
\end{equation*} 
as the boundary of $R$ (see Figure \ref{Fig. Local section in images}(c)), we define $\widehat{R} = R \setminus B$ (see Figure \ref{Fig. Local section in images}(d)). Defining $X_1 = X \cap \widehat{R}$ and $X_2 = X \setminus R$ (see Figure \ref{Fig. Local section in images}(e)-(h)), we obtain a triad $(X,X_1,X_2)$ with the property $X_1 \cap X_2  = \emptyset$. Because $X$, $X_1$, and $X_2$ are subspaces in $\bbR^2$ that are formed by finitely many closed squares in $\bbR^2$, we must have $\cl_X(X_1) \cap \cl_X(X_2) = \emptyset$. 

For example, the local system $(X,X_1,X_2)$ in Figure \ref{Fig. Local section in images} has merging and outer-merging numbers $m_0(X_1;X_2) = 3$, $o_0(X_1;X_2) = 0$, $m_1(X_1;X_2) = 0$, and $o_1(X_1;X_2) = 1$. In~\cite{hu2021sheaf}, we separated a 2D image into disjoint blocks $X_i^{(i)}$ (called a \textit{local patches}) and calculated the local merging numbers $m_0(X_1^{(i)};X_2^{(i)})$ to form a heatmap of the image. In the paper, we mainly focus on the number $o_1(X_1; X_2)$ to approximate the hole positions in a binary image. We show in Section \ref{Sec. Demonstration} how to use the local system described in this section to construct local patches in the image and use them to estimate the 1D holes in the image.

\subsection{Discussion}
\label{Subsec. Discussion}

The organization of the section is as follows. First, we compare the proposed method with possible methods in Section \ref{Subsubsec. Comparison with other methods}. Second, Section \ref{Subsubsec. Size Issues} discusses how to estimate the size and shape of the holes in the topological space through the local system. Finally, Section \ref{Subsubsec. More general systems} discusses local systems composed of n subspaces and their local sections, which will be an important future research direction to promote the theory of this paper.

\begin{figure}
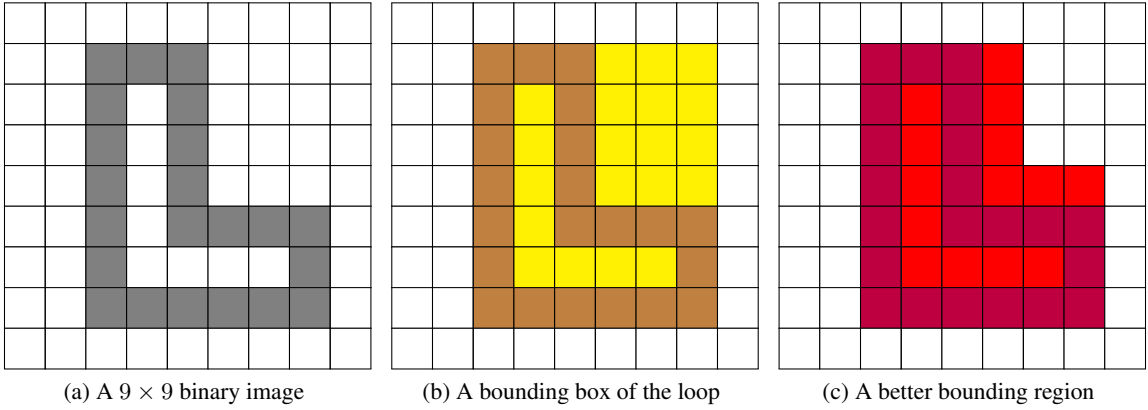
%
  \centering
  \subfloat[A $9 \times 9$ binary image]{
\begin{ytableau}
*(white) & *(white) & *(white) & *(white) & *(white) & *(white) & *(white) & *(white) & *(white)\\
*(white) & *(white) & *(gray) & *(gray) & *(gray) & *(white) & *(white) & *(white) & *(white)\\
*(white) & *(white) & *(gray) & *(white) & *(gray) & *(white) & *(white) & *(white) & *(white)\\
*(white) & *(white) & *(gray) & *(white) & *(gray) & *(white) & *(white) & *(white) & *(white)\\
*(white) & *(white) & *(gray) & *(white) & *(gray) & *(white) & *(white) & *(white) & *(white)\\
*(white) & *(white) & *(gray) & *(white) & *(gray) & *(gray) & *(gray) & *(gray) & *(white)\\
*(white) & *(white) & *(gray) & *(white) & *(white) & *(white) & *(white) & *(gray) & *(white)\\
*(white) & *(white) & *(gray) & *(gray) & *(gray) & *(gray) & *(gray) & *(gray) & *(white)\\
*(white) & *(white) & *(white) & *(white) & *(white) & *(white) & *(white) & *(white) & *(white)\\
  \end{ytableau}}~
  \subfloat[A bounding box of the loop]{
\begin{ytableau}
*(white) & *(white) & *(white) & *(white) & *(white) & *(white) & *(white) & *(white) & *(white)\\
*(white) & *(white) & *(brown) & *(brown) & *(brown) & *(yellow) & *(yellow) & *(yellow) & *(white)\\
*(white) & *(white) & *(brown) & *(yellow) & *(brown) & *(yellow) & *(yellow) & *(yellow) & *(white)\\
*(white) & *(white) & *(brown) & *(yellow) & *(brown) & *(yellow) & *(yellow) & *(yellow) & *(white)\\
*(white) & *(white) & *(brown) & *(yellow) & *(brown) & *(yellow) & *(yellow) & *(yellow) & *(white)\\
*(white) & *(white) & *(brown) & *(yellow) & *(brown) & *(brown) & *(brown) & *(brown) & *(white)\\
*(white) & *(white) & *(brown) & *(yellow) & *(yellow) & *(yellow) & *(yellow) & *(brown) & *(white)\\
*(white) & *(white) & *(brown) & *(brown) & *(brown) & *(brown) & *(brown) & *(brown) & *(white)\\
*(white) & *(white) & *(white) & *(white) & *(white) & *(white) & *(white) & *(white) & *(white)\\
  \end{ytableau}}~
  \subfloat[A better bounding region]{
\begin{ytableau}
*(white) & *(white) & *(white) & *(white) & *(white) & *(white) & *(white) & *(white) & *(white)\\
*(white) & *(white) & *(purple) & *(purple) & *(purple) & *(red) & *(white) & *(white) & *(white)\\
*(white) & *(white) & *(purple) & *(red) & *(purple) & *(red) & *(white) & *(white) & *(white)\\
*(white) & *(white) & *(purple) & *(red) & *(purple) & *(red) & *(white) & *(white) & *(white)\\
*(white) & *(white) & *(purple) & *(red) & *(purple) & *(red) & *(red) & *(red) & *(white)\\
*(white) & *(white) & *(purple) & *(red) & *(purple) & *(purple) & *(purple) & *(purple) & *(white)\\
*(white) & *(white) & *(purple) & *(red) & *(red) & *(red) & *(red) & *(purple) & *(white)\\
*(white) & *(white) & *(purple) & *(purple) & *(purple) & *(purple) & *(purple) & *(purple) & *(white)\\
*(white) & *(white) & *(white) & *(white) & *(white) & *(white) & *(white) & *(white) & *(white)\\
  \end{ytableau}}~
  \caption{An illustration of bounding regions of the hole structure. (a) A binary image that contains a single 1-dimensional loop structure. (b) A rectangular bounding box formed by yellow and brown pixels. (c) A more compact bounding region formed by red and purple pixels.}
  \label{Fig. holes in a binary image}
\end{figure}

\subsubsection{Comparison with other methods}
\label{Subsubsec. Comparison with other methods}

To detect the local regions that contain pores in an $m \times n$ binary image, some naive methods can be used to tackle this question. For example, one can search every subfigure of the given image and compute whether it contains hole structures. However, it is generally infeasible to compute all the
\begin{equation*}
    \binom{m}{2} \cdot \binom{n}{2} = \frac{(m^2-m)(n^2-n)}{4},
\end{equation*}
subfigures for large $m$ and $n$. Except for the computational complexity, covering an irregular hole costs a large bounding rectangle and makes the estimation less precise and compact (see Figure \ref{Fig. holes in a binary image}). 

Following our previous approach~\cite{hu2021sheaf}, we have two strategies to reduce computational complexity through the image's local systems and sheaf information. The first is splitting the image into many pairs $(X_1^{(i)}, X_2^{(i)})$ with disjoint $X_1^{(i)}$s and computing each $\mathcal{P}_q(\mathcal{G}_{i,1})$ on the filtration $\mathcal{G}_{i,1}: \emptyset \subseteq X_1^{(i)} \subseteq X_1^{(i)} \cup X_2^{(i)} \subseteq X$. The second is to use a sliding window technique to cover the entire image and compute short-persistent homology for each local window. The second strategy computes the persistent homology $O(mn)$ times,  and pore locations generated by this strategy are usually more refined than the first one. In the paper, we mainly follow the second strategy and show in Section \ref{Sec. Demonstration} that the proposed barcode and local system framework can detect local holes effectively and with more concise bounding regions than bounding boxes (Figure \ref{Fig. holes in a binary image} (b), (c)).

\subsubsection{Size Issues}
\label{Subsubsec. Size Issues}

As we mentioned in Section \ref{Subsec. Persistent Homology}, homological generators generally lack specific geometric properties, such as the size and stability of pores (Figure \ref{Fig. Opening filtration}), which cannot be detected by traditional homology or persistent homology of sub-level set filtration. In recent years, the shape and size of pore structures have become more and more important topics in bioinformatics and material science~\cite{jiang2018pore,jiao2008modeling,vogel2001quantitative,wu2023effect,wei2023distribution,bassu2023microgel}. Recently, some research has shown the potential and advantage of persistent homology in pore size analysis~\cite{ishihara2023effect,yamauchi2023bin,hu2021two,chung2022multi}. 

As far as the field of image processing is concerned, the combination of persistent homology and mathematical morphology has opened up a new research direction for this field~\cite{garin2019topological,hu2021two,chung2022multi,tymochko2020using}. In particular, our approach in~\cite{chung2022multi} applies morphological opening and closing to measure the spatial information of black and white regions in binary images. Mathematical morphological operations can estimate the sizes of image pores, and most of the current work focuses on the global description of such spatial information, such as the number of pores with a specified morphological size and the average image pore size. However, the location information of pores in images is still limited in present methods. Through the discussion in this section, we will see that localized systems can capture both the location and size of pores, providing a richer pore analysis technique.

\begin{theorem}
\label{Theorem: < 3 theorem}
Let $F$ be a field. Let $(X, X_1, X_2)$ be a local system of topological spaces and $q$ a non-negative integer.  Then, non-zero elements in $H_q(X_1)$ and $H_q(X_2)$ in the persistent homology $\mathcal{P}_q(\mathcal{G}_1): 0 \xrightarrow{} H_q(X_1; F) \xrightarrow{} H_q(X_1 \cup X_2; F) \xrightarrow{} H_q(X; F)$ has birth number $< 3$. 
\end{theorem}
\begin{proof}
Suppose $s_1$ is a non-zero element in $H_q(X_1)$, then the birth number of $s_1$ is $1$. On the other hand, the assumption of $\cl_X(X_1) \cap \cl_X(X_2) = \emptyset$ forces that $H_q(X_1 \cup X_2; F) \simeq H_q(X_1;F) \oplus H_q(X_2;F)$ canonically. Then every non-zero element in $H_q(X_2)$ must have barcode $2$.
\end{proof}

\begin{coro}
\label{Corollary of  < 3 theorem}
Let $X$ be a subspace of the $n$-dimensional Euclidean space $\bbZ^n$. Let  $F$ be a field and $q$ a non-negative integer. For every $c \in Z_q(X;F)$ with $[c] \neq 0$ in $H_q(X;F)$, there is a bounded set $X_1 \subseteq X$ and an $X_2 \subseteq X$ such that $\cl_X(X_1) \cap \cl_X(X_2) = \emptyset$ and $c \in H_q(X_1;F)$. In particular, $[c]$ has a barcode $(1,\star)$ in $\mathcal{P}_q(\mathcal{G}_1)$.   
\end{coro}
\begin{proof}
For a chain $c$ in $S_q(X;F)$, we can write $c = \sum_{i = 1}^n \lambda_i \sigma_i$, where $\lambda_i \in F \setminus \{ 0 \}$ and $\sigma_i: \Delta_q \xrightarrow{} X$ is a continuous for each $i$. Recall that the support of $c$ denoted by $|c|$ is defined as the union of the images of the $\sigma_i$. Because $\Delta_q$ is compact, and $\sigma_i$ is continuous, the support of $c$ is a compact subset of $X$.  In particular, the support of $c$ is closed and bounded. Therefore, we may choose a positive number $r_1 $ such that $|c| \subseteq \mathbf{B}(\mathbf{0},r_1)$. Choose $r_2 > r_1$ and set $X_1 = X \cap \mathbf{B}(\mathbf{0},r_1)$ and $X_2 = X \cap \{ \bfx \in \bbR^n : |\bfx| \geq r_2 \}$, then $\cl_X(X_1) \cap \cl_X(X_2) = \emptyset$ and $c \in Z_q(X_1;F)$. By the proof of Theorem \ref{Theorem: < 3 theorem}, $[c]$ has a barcode $(1,\star)$ in $\mathcal{P}_q(\mathcal{G}_1)$.
\end{proof}

\begin{def.}
For convenience, we use $i_q(X_1;X_2)$ to denote the number of barcodes in $\mathcal{P}_q(\mathcal{G}_1)$.  
\end{def.}

Corollary \ref{Corollary of  < 3 theorem} gives us a way to measure the size of $q$-holes ($q > 0$) for any subspace in $\bbR^n$ by choosing the local system appropriately. More precisely, we can choose a bounded subspace $X_1$ of $X$ that contains the hole and is therefore an approximation of the size of the hole. Also, since the location of $X_1$ is known, it also keeps track of the hole location. We also demonstrate in Section \ref{Sec. Demonstration} an application of Corollary \ref{Corollary of  < 3 theorem} to detect the "largest" holes in images.


\subsubsection{More general systems}
\label{Subsubsec. More general systems}

In the paper, we focus on a local system consisting of topological spaces $X, X_1, X_2$ that satisfy $X_1 \subseteq X$, $X_2 \subseteq X$, and $\cl_X(X_1) \cap \cl_X(X_2) = \emptyset$. This system induces a sheaf structure as in Theorem \ref{Theorem: Global section space of the simplest sheaf}, and its global section space can be computed by the persistent homology of the short filtration. Actually, one can consider a more general case consisting of $n+1$ spaces $X, X_1, ..., X_n$ with $\cl_X(X_i) \cap \cl_X(X_j) = \emptyset$ for $i \neq j$. We synthesize the above settings into the following definition and theorem.
\begin{def.}
Let $X$ be a topological space and $X_1, ..., X_n$ by subspaces of $X$ that satisfy $\cl_X(X_i) \cap \cl_X(X_j) = \emptyset$ for $i \neq j$. The $(n+1)$-tuple $(X,X_1, ..., X_n)$ is called a \textbf{local }$\textbf{n}$\textbf{-system} (or an \textbf{admissible} $\textbf{(n+1)}$\textbf{-tuple}) of topological spaces. 
\end{def.}

Furthermore, for an $n$-system $(X_1, ..., X_n)$, we can consider the diagram
\begin{equation*}
\xymatrix@+1.5em{
				& H_q(X_1)
				\ar[dr]^{\rho_{1}}
                &  
                \\
        	& H_q(X_k) \ar@{.}[u]^{} \ar[r]_{\rho_{k}}
				& H_q(X) 
				\\
        	& H_q(X_n) \ar@{.}[u]^{} \ar[ur]_{\rho_{n}}
				& 
}
\end{equation*}
with homologies and induced homomorphisms~\cite{hu2020brief}. As above, we define its global section space by 
\begin{equation*}
    \Gamma = \left\{ (s_1, s_2, ..., s_n) \in \prod_{i=1}^n H_q(X_i) : \rho_i(s_i) = \rho_j(s_j) \text{ for } i, j \in \{ 1,2, ..., n\} \right\}.
\end{equation*}

\begin{theorem}
\label{Theoerm: Global section thm of n-system}
Let $(X,X_1, ..., X_n)$ be a local $n$-system of topological spaces. Let $R$ be a commutative ring with identity and $q \geq 0$ a non-negative integer. Consider the sheaf structure
\begin{equation*}
\xymatrix@+1.5em{
				& H_q(X_1;R)
				\ar[dr]^{\rho_{1}}
                &  
                \\
        	& H_q(X_k;R) \ar@{.}[u]^{} \ar[r]_{\rho_{k}}
				& H_q(X;R) 
				\\
        	& H_q(X_n;R) \ar@{.}[u]^{} \ar[ur]_{\rho_{n}}
				& 
}
\end{equation*}
of $R$-modules and module homomorphisms and the homomorphism
\begin{equation}
\label{Eq. Global section homomorphism}
\bigoplus_{i = 1}^n H_{q}(X_i;R) \xrightarrow{ \ \ \phi \ \ } \bigoplus_{i = 2}^n H_{q}(X;R), \ \ (s_i)_{i = 1}^n \longmapsto (\rho_1(s_1) - \rho_i(s_i))_{i = 2}^n.
\end{equation}
Let $\Gamma$ be the global section space of the sheaf. Then $\Gamma$ is the kernel of $\phi$. 
\end{theorem}
\begin{proof}
An $n$-tuple $(s_i)_{i = 1}^n$ is a global section if and only if $\rho_i(s_i) - \rho_j(s_j) = 0$ for every $i, j \in \{ 1,2, ..., n \}$. If it holds, then $\rho_1(s_1) - \rho_j(s_j) = 0$ for every $j \in \{ 2, ..., n \}$.  Conversely, suppose $\rho_1(s_1) - \rho_j(s_j) = 0$ for every $j \in \{ 2, ..., n \}$, then $\rho_i(s_i) - \rho_j(s_j) = \rho_i(s_i) - \rho_1(s_1) + \rho_1(s_1) - \rho_j(s_j) = 0$ for every $i, j \in \{ 2, ..., n \}$, as desired.
\end{proof}

Although Theorem \ref{Theoerm: Global section thm of n-system} provides a way to sculpt the global section space, it still has a limitation in computation. When $n = 2$, and $R = F$ is a field, the homomorphism in~\eqref{Eq. Global section homomorphism} and $H_q(X_1 \cup X_2;F) \xrightarrow{} H_q(X;F)$ have the same image.  In this case, the approximation developed in Theorem \ref{Theorem: Local section approximation} and Theorem \ref{Theroem for barcode (3,inf)} is available. However, the map in~\eqref{Eq. Global section homomorphism} and the canonical one $H_q(X_1 \cup X_2;F) \xrightarrow{} H_q(X;F)$ are not coincident and can not be calculated by counting the barcodes in the short filtration. Computing global sections through persistence barcodes is one of our future research directions.

\begin{figure*}
    \begin{center}
    \subfloat[Input image]{\includegraphics[width=40mm]{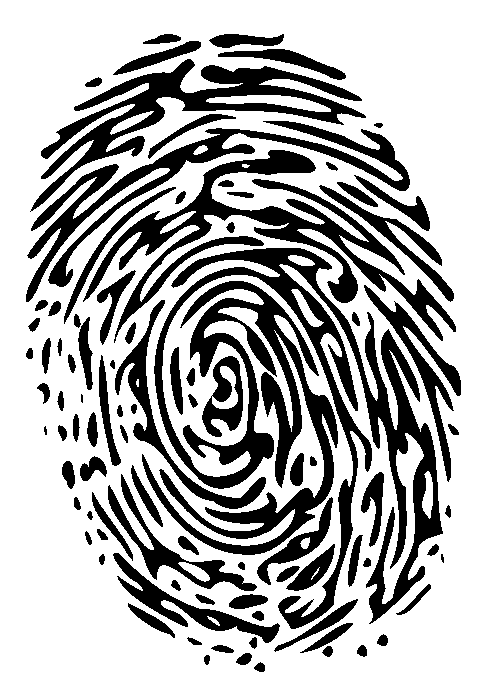}}~
    \subfloat[Marked holes]{\includegraphics[width=40mm]{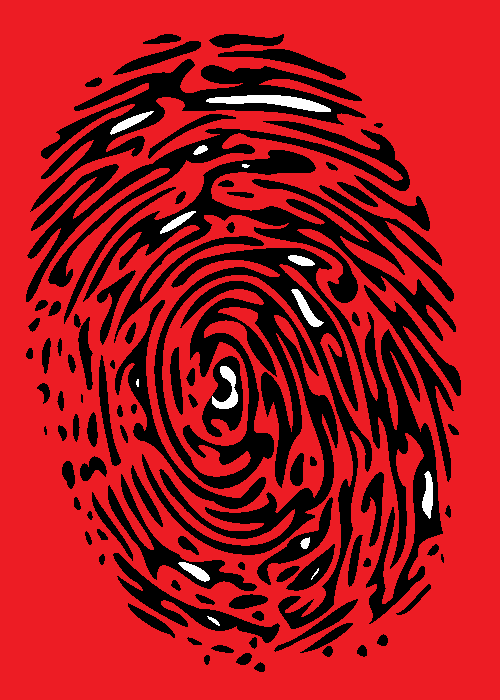}}~
    \subfloat[Output heatmap]{\includegraphics[width=40mm]{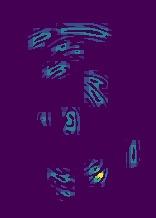}}~
    \subfloat[Location estimation]{\includegraphics[width=40mm]{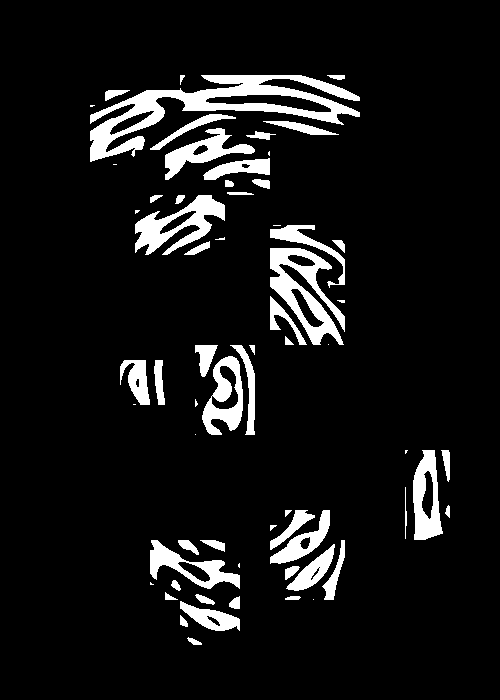}}
    \end{center}
    \caption{A demonstration of Algorithm \ref{alg:  hole structure detection algorithm} on a binary image. (a) An illustration of a  $500 \times 700$ fingerprint image, where the Betti pair of the image is $(\beta_0,\beta_1) = (92,14)$. (b) The marked white regions as the $14$ holes of the fingerprint image.  (c) The output heatmap by Algorithm \ref{alg:  hole structure detection algorithm}. (d) The non-zero parts of the output heatmap, form an estimation for the hole positions. Here we choose the window $R$ as a $30 \times 30$ square with step $k = 15$.}
    \label{Fig. Demostration_1}
 \end{figure*}

\section{Demonstration on Digital Images}
\label{Sec. Demonstration}

Hole structures in images can be subtle and complicated. As shown in Figure \ref{Fig. Demostration_1}(a)\footnote{Karen Arnold has released this “Fingerprint Clipart” image under a Public Domain license. \url{https://www.publicdomainpictures.net/en/view-image.php?image=462168&picture=fingerprint-clipart}}, although many white areas appear in the image as porosity structures or closed voids, many of these areas connect to the white background and thus are not actual holes. In order to detect the image's hole positions, we propose Algorithm  \ref{alg:  hole structure detection algorithm} to approximate the holes' geometric locations using the above theory and sliding window technique. Figure \ref{Fig. Demostration_1} is a demonstration of Algorithm  \ref{alg:  hole structure detection algorithm} on a binary image. We can see that all the holes in the image are detected by the output heatmap as Figure~\ref{Fig. Demostration_1}(c). 

We note that Line \ref{line 6 in Algorithm 1} in Algorithm \ref{alg:  hole structure detection algorithm} considers both $i_1(X_1;X_2)$ and $o_1(X_1;X_2)$. If $i_1(X_1; X_2) \neq 0$, it means that $X_1$ contains some holes in $X$ and already a bounding box of certain holes (Corollary \ref{Corollary of  < 3 theorem}). On the other hand, $i_1(X_1; X_2)$ records whether (part of) the black pixels in $X_1$ can contribute to hole structures and the number of these structures. Therefore, the sum of $i_1(X_1;X_2)$ and $o_1(X_1;X_2)$ can estimate whether $X_1$ nears a hole structure in $X$.

\begin{algorithm}[H]
\renewcommand{\algorithmicrequire}{\textbf{Input:}}
\renewcommand{\algorithmicensure}{\textbf{Output:}}
\begin{algorithmic}[1]
\Require Binary image $f: P \xrightarrow{} \{0,1\}$ on a rectangle $R$, $X = f^{-1}(0)$, an $n \times n$ square window $R$, and a sliding step $k$.
\Ensure A function $H: P \rightarrow \bbR$ as a heatmap of $f$. The heatmap estimates the hole locations in image $f$. A point in $P$ with a high heat value is more possible as a part of a hole.
\State Denote $P = ([0,a] \times [0,b]) \cap \bbZ^2$ and $R = ([0,n] \times [0,n]) \cap \bbZ^2$. Define $B$ and $\widehat{R}$ as in Section \ref{Subsec. Local Systems in Binary Images}.
\State Set $H: P \rightarrow \bbR$ as the zero function.
\For{$i \in \{ 0,1,..., a\}$ and $j \in \{ 0,1,..., b\}$}
\If{$(i \cdot k, j \cdot k) + R \subseteq P$}
    \State Set $X_1 = ((i,j) + \widehat{R}) \cap X$ and $X_2 = X \setminus ((i,j) + R)$
    \State Compute $\mathcal{M} = i_1(X_1;X_2) + o_1(X_1;X_2)$
    \label{line 6 in Algorithm 1}
    \State Define $H': P \xrightarrow{} \{ 0,1 \}$ as follows:
    \begin{equation*}
      H'(\bfx) = \begin{cases}
 	H(\bfx) + \mathcal{M} &\null\text{ if } \bfx \in (i,j) + \widehat{R}, \\
 	0  &\null\text{ otherwise.}
 	\end{cases}  
    \end{equation*}
    \State $H \xleftarrow{} H'$
\Else
    \State \textbf{continue}
\EndIf
\EndFor \\
\Return $H' \cdot (1 - f)$
\end{algorithmic}          
    \caption{The hole structure detection algorithm.} 
    \label{alg:  hole structure detection algorithm}
\end{algorithm}

Apart from the task of detecting the location of holes in an image, recognizing the size or shape of holes is also an interesting one. As introduced in Section \ref{Subsubsec. Size Issues}, local windows that contain holes will produce $(1,2)$ barcodes in the corresponding short persistent homology. Based on the observation, we can modify Algorithm \ref{alg:  hole structure detection algorithm} by considering the information on the size of local windows and changing the $\mathcal{M}$ value to tackle this task. As follows, we propose Algorithm \ref{alg:  hole structure detection algorithm-2} as a modification of Algorithm \ref{alg:  hole structure detection algorithm}. We note here that the we implement these two algorithms in {\tt{Python}} with the {\tt{Gudhi}} package~\cite{gudhi:urm}.

\begin{algorithm}[H]
\renewcommand{\algorithmicrequire}{\textbf{Input:}}
\renewcommand{\algorithmicensure}{\textbf{Output:}}
\begin{algorithmic}[1]
\Require Binary image $f: P \xrightarrow{} \{0,1\}$ on a rectangle $R$, $X = f^{-1}(0)$, an $n \times n$ square window $R$, and a sliding step $k$.
\Ensure A function $H: P \rightarrow \bbR$ as a heatmap of $f$. The heatmap estimates the hole locations in image $f$. A point in $P$ with a high heat value is more possible as a part of a ``large'' hole.
\State Denote $P = ([0,a] \times [0,b]) \cap \bbZ^2$ and $R = ([0,n] \times [0,n]) \cap \bbZ^2$. Define $B$ and $\widehat{R}$ as in Section \ref{Subsec. Local Systems in Binary Images}.
\State Set $H: P \rightarrow \bbR$ as the zero function.
\For{$i \in \{ 0,1,..., a\}$ and $j \in \{ 0,1,..., b\}$}
\If{$(i \cdot k, j \cdot k) + R \subseteq P$}
    \State Set $X_1 = ((i,j) + \widehat{R}) \cap X$ and $X_2 = X \setminus ((i,j) + R)$
    \State Compute $\mathcal{M} = {\rm vol}(R) \cdot (o_1(X_1;X_2) - i_1(X_1;X_2))$  \Comment{The only difference to Algorithm \ref{alg:  hole structure detection algorithm}}
    \State Define $H': P \xrightarrow{} \{ 0,1 \}$ as follows:
    \begin{equation*}
      H'(\bfx) = \begin{cases}
 	H(\bfx) + \mathcal{M} &\null\text{ if } \bfx \in (i,j) + \widehat{R}, \\
 	0  &\null\text{ otherwise.}
 	\end{cases}  
    \end{equation*}
    \State $H \xleftarrow{} H'$
\Else
    \State \textbf{continue}
\EndIf
\EndFor \\
\Return $H' \cdot (1 - f)$
\end{algorithmic}          
    \caption{The hole size estimation algorithm.} 
    \label{alg:  hole structure detection algorithm-2}
\end{algorithm}

We note that the only difference between Algorithm \ref{alg:  hole structure detection algorithm} and Algorithm \ref{alg:  hole structure detection algorithm-2} is the value of $\mathcal{M}$ in line 6 of both algorithms. As we discussed above, $i_1(X_1; X_2) > 0$ means that some holes are bounded by the area $X_1$. Due to this, we apply $-{\rm vol}(R) \cdot i_1(X_1; X_2)$ as the punishment term of the region's heat value. In addition, since punishment will produce negative values, if a certain area has many fine holes, the heat function will have a high negative value in this area and also estimates the location of these holes.

\begin{figure*}
    \begin{center}
    \subfloat[Input image]{\includegraphics[width=26mm]{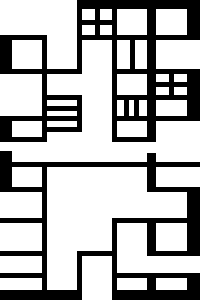}}~
    \subfloat[$50 \times 50$]{\includegraphics[width=26mm]{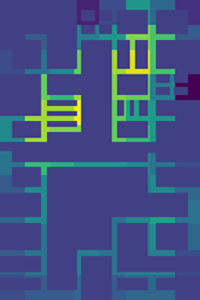}}~
    \subfloat[$100 \times 100$]{\includegraphics[width=26mm]{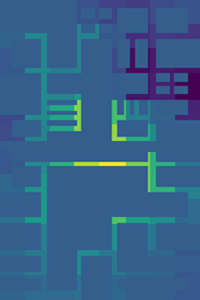}}~
    \subfloat[$150 \times 150$]{\includegraphics[width=26mm]{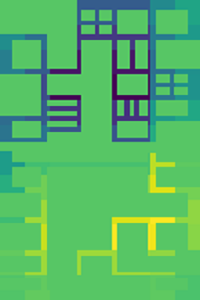}}~
    \subfloat[Sum of (b)-(d)]{\includegraphics[width=26mm]{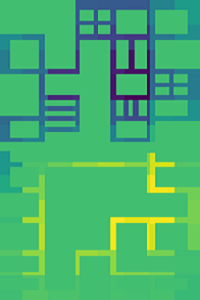}}~
    \subfloat[Estimation]{\includegraphics[width=26mm]{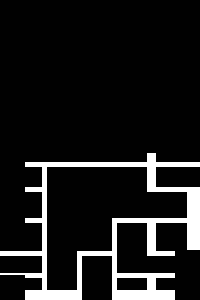}}
    \end{center}
    \caption{A demonstration of Algorithm \ref{alg:  hole structure detection algorithm-2} on a binary image. (a) An illustration of a $200 \times 300$ binary image with Betti pair $(\beta_0,\beta_1) = (2,28)$. (b), (c), and (d) are the output heatmaps of Algorithm \ref{alg:  hole structure detection algorithm-2} with local windows of sizes $50 \times 50$, $100 \times 100$, and $150 \times 150$, respectively. (e) The sum heatmap of (b), (c), and (d). (e) The non-zero parts of the output heatmap, form an estimation for the hole positions. Here we choose the step $k = 25$.}
    \label{Fig. Demostration_2}
 \end{figure*}

Figure \ref{Fig. Demostration_2} demonstrates Algorithm  \ref{alg:  hole structure detection algorithm-2} on a binary image with hole structures in different sizes. Within different local windows ($50 \times 50$, $100 \times 100$, and $150 \times 150$ squares), Algorithm  \ref{alg:  hole structure detection algorithm-2} gives different attention to these holes. We notice that small holes may get more attention from Algorithm  \ref{alg:  hole structure detection algorithm-2} (Figure \ref{Fig. Demostration_2}(b)) since there are many holes next to each other, and hence the shared edges will get a higher heat value. Keeping enlarging the size of the local window, we observe that smaller holes in the image would get more punishment in heat values. The sum of the three heatmaps summarizes the ``importances'' of black pixels in the image. Finally, we see that Algorithm \ref{alg:  hole structure detection algorithm-2} successfully approximates the position of the ``largest hole'' by the thresholding method.

However, real pore structure data may be more complex than the images presented in the paper. Methods to use barcode information in pore structure analysis, such as the choice of the penalty function, still need research and development, which is our main development work in the future.

\section{Conclusion}
\label{Sec. Conclusion} 
To summarize the paper, we propose a local system and local persistent homology framework to study the local merging relations in an arbitrary topological space. By using the merging and out-merging numbers of local regions, we propose an algorithm to detect the sizes and positions of pores in the image. Although the demonstration focuses on digital images, the framework can be adapted to any topological space with local systems. We also look forward to applying this framework to point-cloud data, especially its applications in crystalline data analysis.

\section*{Acknowledgement}
Most of the work in this article was completed by the author during his doctoral study at National Taiwan Normal University (2016-2022). The author would like to thank Dr. Chun-Chi Lin (NTNU) and Dr. Yu-Min Chung (Eli Lilly and Company), the author's doctoral supervisors, for their comments and suggestions on the work. Especially, Dr. Yu-Min Chung provided many suggestions for studying the geometric meaning of persistent barcodes and the global sections of the local system, making the discussion more fruitful and rigorous. The author would also like to thank Dr. Kelin Xia, the author's postdoctoral supervisor at Nanyang Technological University. The author got a lot of inspiration from the discussion with Dr. Xia so that this paper can have more research directions, such as a more detailed study of the geometry of cellular sheaves, and more possible applications.

\bibliographystyle{abbrv}  
\bibliography{references}

\end{document}